\documentclass[12pt,a4paper]{amsart}

\usepackage{amsfonts}
\usepackage{amsmath}
\usepackage{amssymb}
\usepackage{amsthm}
\usepackage{mathrsfs}
\usepackage{color}
\usepackage{bbm}
\usepackage{MnSymbol}
\usepackage{extarrows}

\usepackage{graphicx}
\usepackage[all]{xy}

\usepackage{enumerate}

\usepackage{hyperref}

\DeclareMathAlphabet{\mathpzc}{OT1}{pzc}{m}{it}

\makeatletter
\@namedef{subjclassname@2020}{\textup{2020} Mathematics Subject Classification}
\makeatother

\newcommand{\f}{\varphi}

\newcommand{\aA}{{\mathcal{A}}}
\newcommand{\bB}{{\mathcal{B}}}
\newcommand{\cC}{{\mathcal{C}}}
\newcommand{\dD}{{\mathcal{D}}}

\newcommand{\fF}{{\mathcal{F}}}

\newcommand{\hH}{{\mathcal{H}}}

\newcommand{\nN}{{\mathcal{N}}}
\newcommand{\lL}{{\mathcal{L}}}
\newcommand{\oO}{{\mathcal{O}}}
\newcommand{\pP}{{\mathcal{P}}}

\newcommand{\sS}{{\mathcal{S}}}

\newcommand{\wW}{{\mathcal{W}}}

\newcommand{\zZ}{{\mathcal{Z}}}

\newcommand{\kk}{\mathbbm{k}}

\newcommand{\la}{\lambda}

\newcommand{\tr}{{\mbox{$t$-struc}\-tu\-r}}

\newcommand{\bcdot}{{\mbox{\boldmath{$\cdot$}}}}
\newcommand{\opp}{\mathop{\textrm{op}}\nolimits}	
\newcommand{\Hom}{\mathop{\textrm{Hom}}\nolimits}
\newcommand{\coker}{\mathop{\textrm{coker}\,}\nolimits}
\newcommand{\Def}{\mathop{\textrm{Def}}\nolimits}
\newcommand{\ncDef}{\mathop{\textrm{ncDef}}\nolimits}
\newcommand{\Ext}{\mathop{\textrm{Ext}}\nolimits}
\newcommand{\End}{\mathop{\textrm{End}}\nolimits}

\newcommand{\Coh}{\mathop{\textrm{Coh}}\nolimits}
\newcommand{\QCoh}{\mathop{\textrm{QCoh}}\nolimits}

\newcommand{\Rex}{\mathop{\textrm{Rex}}\nolimits}

\newcommand{\Spec}{\mathop{\textrm{Spec}}\nolimits}

\newcommand{\modu}{\mathop{\textrm{mod--}}\nolimits}
\newcommand{\Modu}{\mathop{\textrm{Mod--}}\nolimits}

\newcommand{\Per}[1]{\mathop{{}^{#1}\textrm{Per}}\nolimits}
\newcommand{\Id}{\mathop{\textrm{Id}}\nolimits}

\newcommand{\Fin}{\mathop{\textrm{Fin}}\nolimits}
\newcommand{\Df}{\mathop{\textrm{Df}}\nolimits}

\newcommand{\Bas}{\mathop{\textrm{Bas}}\nolimits}
\newcommand{\Ind}{\mathop{\textrm{Ind}}\nolimits}

\newcommand{\Sets}{\mathop{\textrm{Sets}}\nolimits}

\newcommand{\rad}{\mathop{\mathpzc{R}}\nolimits}
\newcommand{\ngbh}{\mathop{\mathpzc{Ex}}\nolimits}

\newcommand{\wt}[1]{\widetilde{#1}}
\newcommand{\ol}[1]{\overline{#1}}

\newtheorem*{THM*}{Theorem}
\newtheorem*{COR*}{Corollary}
\newtheorem*{PROP*}{Proposition}
\newtheorem*{LEM*}{Lemma}

\newtheorem{LEM}{Lemma}[section]

\newtheorem{THM}[LEM]{Theorem}
\newtheorem{PROP}[LEM]{Proposition}
\newtheorem{COR}[LEM]{Corollary}

\theoremstyle{definition}
\newtheorem{EXM}[LEM]{Example}
\newtheorem{REM}[LEM]{Remark}
\newtheorem{DEF}[LEM]{Definition}

\begin{document}
	\title{Categorifying non-commutative deformations}
	\author{Agnieszka Bodzenta}
	\author{Alexey Bondal}
	\date{\today}
	\address{Agnieszka Bodzenta\\
		Institute of Mathematics, 
		University of Warsaw \\ Banacha 2 \\ Warsaw 02-097,
		Poland} \email{A.Bodzenta@mimuw.edu.pl}
		
	\address{Alexey Bondal\\
		Steklov Mathematical Institute of Russian Academy of Sciences, Moscow, Russia, and \newline
		AGHA Laboratory, Moscow Institute of Physics and Technology, Russia, and\newline
		Kavli Institute for the Physics and Mathematics of the Universe (WPI), The University of Tokyo, Kashiwa, Chiba 277-8583, Japan} 
	\email{bondal@mi-ras.ru}
	
	\begin{abstract}
		We define the functor $\textrm{ncDef}_{(Z_1,\ldots,Z_n)}$ of non-commutative deformations of an $n$-tuple of objects in an arbitrary $\mathbbm{k}$-linear abelian category $\mathcal{Z}$. 
	In our categorified approach, we view the underlying spaces of infinitesimal flat deformations
	as Deligne finite categories, i.e.  finite length abelian categories admitting projective generators, with $n$ isomorphism classes of simple objects.

	More generally, we define the functor $\textrm{ncDef}_{\zeta}$ of non-commutative deformations of an exact functor $\zeta \colon \mathcal{A} \to \mathcal{Z}$ of abelian categories. 
	Here the role of an infinitesimal non-commutative thickening of $\aA$
	is played by an abelian category $\bB$ containing $\aA$ and such that $\aA$ generates $\bB$ by extensions. The functor $\ncDef_\zeta$ assigns to such $\bB$ the set of equivalence classes of exact functors $\bB \to \zZ$ which extend $\zeta$.
	We prove that an exact  functor on an infinitesimal extension is fully faithful  if and only if it is fully faithful on the first infinitesimal neighbourhood.
	
	We show that if $\zeta$ is fully faithful, then the functor $\textrm{ncDef}_{\zeta}$ is ind-represented by the extension closure of the essential image of $\zeta$.
		
	We prove that for a flopping contraction $f\colon X\to Y$ with the fiber over a closed point $C = \bigcup_{i=1}^n C_i$, where $C_i$'s are irreducible curves, $\{\mathcal{O}_{C_i}(-1)\}$   is the set of simple objects in the null-category for $f$. We conclude that the null-category
		ind-represents the functor $\textrm{ncDef}_{(\mathcal{O}_{C_1}(-1),\ldots,\mathcal{O}_{C_n}(-1))}$. 
	\end{abstract}
	
	\maketitle
	
	\tableofcontents

\section{Introduction}

    In this paper, we develop a categorical point of view on non-commutative deformation theory.
	The classical commutative infinitesimal deformation theory of a mathematical object defined over a base field $\kk$ is governed by the so-called {\em deformation functor} on the  category of commutative local Artinian algebras  with values in the category of sets \cite{Schle}. The functor assigns to such an algebra $A$ all possible flat deformations of the object over $\Spec A$.
	O.~A.~Laudal initiated the study of non-commutative infinitesimal deformations of a module over an algebra by replacing Artinian algebras by non-commutative Artinian algebras \cite{Lau}. 
	
	When deforming a collection of objects, a new phenomenon is observed: in contrast to the commutative case, a non-commutative deformation space can be made more complicated than just the product of the deformation spaces for the objects deformed separately. In a sense, objects interact in the deformation like entangled particles in Quantum Mechanics. This interaction of objects is governed by their pairwise $\Ext$ groups. For non-commutative deformations of $n$ modules over a $\kk$-algebra, O.~A.~Laudal suggested to replace Artinian algebras by $n$-pointed algebras \cite{Lau}, which are basic algebras with a fixed splitting $\kk^{\oplus n}\to A$ of the maximal semi-simple quotient $A/\textrm{rad A} \simeq \kk^{\oplus n}$.

Our proposal is to replace the category of $n$-pointed algebras by the category $\Df_n$ of Deligne finite categories, by which we mean finite length abelian categories with $n$ simple objects, $\{ S_i\}_{i=1}^n$, satisfying  the finiteness condition on dimension of Hom and $\Ext^1$ groups and possessing a projective generator (\emph{cf.} \cite{Del}). 
Also we consider a bigger category $\wt{\Df_n}$, whose objects are similar categories with the condition on the existence of a projective generator omitted. These categories play the role of ind-representing objects for deformation functors, analogous to complete local algebras that pro-represent deformation functors in the classical deformation theory of Schlessinger \cite{Schle}.

Our approach is based on the observation that, for a commutative Artinian algebra $A$ and a scheme $X$, a deformation along $\Spec A$ of a quasi-coherent sheaf on $X$, i.e. a sheaf on $X\times \Spec A$ flat over $\Spec A$, is uniquely determined by an \emph{exact} functor $\modu A\to \QCoh(X)$, where $\modu A$ is the category of finite dimensional  right $A$-modules. This supports the idea to define a non-commutative infinitesimal deformation of a finite collection of $n$ objects $\{ Z_i\}_{i=1}^n$ in an abelian category $\mathcal{Z}$ as an exact functor $\aA \to \mathcal{Z}$ out of a Deligne finite category $\aA$, which, for any $i$, takes the  simple object $S_i$ in $\aA$ to $Z_i$. 

The \emph{non-commutative deformation functor} of an  ordered collection $(Z_1,\ldots,Z_n)$ of objects in $\mathcal{Z}$
$$
\ncDef_{(Z_1,\ldots,Z_n)}: \Df_n^{\opp} \to \Sets
$$
assigns to $\aA$ the set of isomorphism classes of these exact functors. A further extension of the theory where the deformation functor is a 2-functor out of the \emph{2-category} $\Df_n$ to the 2-category of groupoids will be discussed in a subsequent paper. In particular, it will be explained that the 2-categorical abbreviation of the $\infty$-category deformation theory for collections of objects in an abelian category does not lose any information of the latter. 

In the categorical approach of this paper, there is no need of the choice of the splitting for the semi-simple part of the algebra, as it was for $n$ pointed algebras. Instead, we consider the category ${\rm Bas}_n$ of basic unital finite dimensional $\kk$ algebras such that $A/\textrm{rad}A \simeq \kk^n$ together with a full order on simple modules, i.e. a homomorphism $q_A \colon A \to \kk^n$ which induces an isomorphism $A/\textrm{rad} A \simeq \kk^n$.
Morphism in $\textrm{Bas}_n$ are conjugacy classes of homomorphisms $\f \colon A \to A'$ satisfying $q_{A'} \circ \f = q_A$. Here, we say that $\f, \psi \colon A \to A'$ are conjugate if there exists an invertible $u \in A'$ with $q_{A'}(u) =1$, such that $\f = u \psi u^{-1}$. Note that morphisms in this category are not just homomorphisms of algebras. 

According to Theorem \ref{thm_equiv_of_cat},  the category ${\rm Bas}_n$ is equivalent to the  category opposite to $\Df_n$. The equivalence takes algebra $A$ from ${\rm Bas}_n$ to the category $\modu A$.

It turns out that the categorified deformation theory of this paper is good for constructing the `universal deformation', which in categorical terms means the  (ind)representability of the non-commutative deformation functor. Namely, we prove the ind-representability under the condition that the collection of objects that we deform is simple. Following Y. Kawamata \cite{Kaw5}, we say that a collection $\{Z_i\}_{i=1}^n$ of objects in a $\kk$-linear abelian category $\zZ$ is \emph{simple}, if
$$
\dim_{\kk} \Hom_{\zZ}(Z_i,Z_j) = \delta_{ij}.
$$
We denote by $\fF(\{Z_i\})$ the smallest full  subcategory in $\zZ$ closed under extensions and containing $\{Z_1,\ldots, Z_n\}$. Then we prove

\vspace{0.3cm}
\noindent
\textbf{Theorem A.}
\textit{
	(see Theorem \ref{thm_pro-represent_for_simple_coll})
	Let $(Z_1,\ldots, Z_n)$ be an ordered  simple collection of objects in abelian category $\zZ$. If $\dim_{\kk} \Ext^1(\bigoplus Z_i, \bigoplus Z_i)$ is finite, then $\fF(\{Z_i\})$ is the object in $\wt{\Df_n}$ 
that ind-represents the non-commutative deformation functor $\ncDef_{(Z_1,\ldots,Z_n)}$.}
\vspace{0.3cm}

Note that  the (pro)representability of the deformation functor does not hold for the category of $n$-pointed algebras, mostly because morphisms in the underlying category were homomorphisms of algebras and not the conjugacy classes of those. Instead, there exists the hull, i.e. a surjective morphism from a representable functor to the relevant deformation functor \cite{Lau}.

The functor $\ncDef_{(Z_1,\ldots,Z_n)} \colon \Df_n^{\opp} \to \Sets$ can be viewed as a deformation of an (exact) functor $T_{(Z_1,\ldots,Z_n)} \colon \modu \kk^n \to \zZ$, i.e. of  an additive functor which maps the $i$'th simple $\kk^n$-module to $Z_i$. Now we replace $\modu \kk^n$ by an  arbitrary abelian category $\aA$. We introduce the category $\Fin_{\aA}$ of abelian categories which are infinitesimal extensions of $\aA$, i.e. categories $\bB$ containing $\aA$ as a full subcategory and such that the extension closure of $\aA$ is $\bB$.

For an exact functor $\zeta \colon \aA \to \zZ$, we consider  a non-commutative deformation functor $\ncDef_\zeta\colon \Fin_{\aA}^{\opp} \to \Sets$. It assigns to $\bB \in \Fin_{\aA}$ the set of equivalences classes of exact functors $F \colon \bB \to \zZ$ together with an isomorphism $F|_{\aA} \xrightarrow{\simeq} \zeta$.

\vspace{0.3cm}
\noindent
\textbf{Theorem B.}
\textit{
	(see Theorem \ref{thm_ind_rep_fully_faith})
	Let $\zeta \colon \aA \to \zZ$ be a fully faithful functor. Then the extension closure $\ngbh^\zZ(\zeta(\aA))$ of the full subcategory $\zeta(\aA) \subset \zZ$ ind-represents the non-commutative deformation functor $\ncDef_{\zeta}$.}
\vspace{0.3cm}

We give a criterion of the fully faithfulness of an exact functor $\Phi \colon \aA\to \zZ$ provided $\aA$ is the extension closure of a full subcategory $\aA_0 \subset \aA$.  We show that $\Phi$ is fully faithful if and only if its restriction to the first infinitesimal neighbourhood of $\aA_0$ is, see Theorem \ref{thm_ff_criterion}.

The original geometric motivation for this paper was to better understand the relation between flops and non-commutative deformations. 
For the case of a flopping contraction of a $(0,-2)$-curve $C$ on a threefold $X$, Y. Toda  proposed in \cite{Toda3} an interpretation of the flop-flop functor as the spherical twist with respect to a spherical functor $\dD^b({\rm mod-}A)\to \dD^b(X)$, where $A=\kk[x]/x^n$. Also he gave an interpretation for ${\rm Spec}A$ as the moduli of (commutative) deformations of $\oO_C(-1)$ in $\Coh (X)$. 

This idea was further developed by W. Donovan and M. Wemyss in \cite{DonWem}. They interpreted the flop-flop functor for a more general threefold flopping contraction with an irreducible flopping curve $C$ in terms of the so-called contraction algebra, which is a hull for the non-commutative deformation functor for $\oO_C(-1)$ in the sense of Laudal.  
Using M.~Van den Bergh's non-commutative resolution which gives an equivalence $\dD^b(X) \xrightarrow{\simeq} \dD^b(\Lambda)$, for a suitable non-commutative algebra $\Lambda$, W.~Donovan and M.~Wemyss transferred their analysis into the pure algebraic framework of $\Lambda$ representations. They considered the contraction algebra $\Lambda_{con}$, which is a quotient of $\Lambda$ by an ideal $I$, and proved that $I$ regarded as a $\Lambda - \Lambda$ bimodule defines a  self-equivalence of 
$\dD^b(\Lambda)$. When this equivalence is transported into $\dD^b(X)$ it coincides with the flop-flop functor. In subsequent paper \cite{DonWem1}, the authors considered the case of reducible flopping curve $C$ on a threefold $X$ and described relations between auto-equivalences for flops with centers in different irreducible components of $C$.

Consider a more general flopping contraction $f\colon X\to Y$ of Gorenstein varieties with fibers of dimension bounded by 1. Assume that variety $Y$ has canonical hypersurface singularities of multiplicity 2 (\emph{cf.} \cite{VdB}). In \cite{BodBon}, we introduced the \emph{null-category}
\begin{equation*}
\mathscr{A}_f =\{E\in \Coh(X)\,|\, Rf_*(E) = 0\}
\end{equation*} 
and proved that under the above assumptions the `flop-flop' equivalence is the \emph{spherical twist} for the spherical functor $\dD^b(\mathscr{A}_f) \to \dD^b(\Coh(X))$, the derived functor of the embedding $\mathscr{A}_f\to \Coh(X)$.

In this paper, we relate the null-category to non-commutative deformations. 
Theorem \ref{thm_simple_objects_in_A_f} states that simple objects in $\mathscr{A}_{f}$ are $\oO_{C_{y,i}}(-1)$, where $C_{y,i} \simeq \mathbb{P}^1$ run over irreducible components  of 1-dimensional fibers $C_y$ for $f$ over all closed points of $y$. Fix one closed point $y \in Y$ with 1-dimensional fiber $C_y$ over it and consider the \emph{null-category with support in the fiber} 
	\begin{equation*}
	\mathscr{A}_{f,C_y}= \{E\in \mathscr{A}_f \,|\, \textrm{Supp }E \subset C_y\}
	\end{equation*}

\vspace{0.3cm}
\noindent
\textbf{Theorem C.} \textit{ (see Theorem \ref{thm_universal_def})
The null-category $\mathscr{A}_{f,C_y}$ is an object of $\wt{\Df}_n$ that ind-represents the non-commutative deformation functor of the collection $\{\oO_{C_{y,i}}(-1)\}$.}
\vspace{0.3cm}

If we assume further that $\dim X = \dim Y =3$ and $Y$ is the spectrum of a noteherian local ring, then any object of $\mathscr{A}_f$ is supported on the fiber $C$ of $f$ over the unique closed point of $Y$, i.e. $\mathscr{A}_{f} \simeq \mathscr{A}_{f,C}$. Moreover, under these assumptions the category $\mathscr{A}_f$ is Deligne finite, hence we have

\vspace{0.3cm}
\noindent
\textbf{Corollary}
\textit{
(see Corollary \ref{cor_univ_def_dim_3})
Let $R$ be a noetherian local ring of dimension 3. If $f\colon X\to \Spec R$ is a flopping contraction, then 
the category $\mathscr{A}_{f}$ 
represents the non-commutative deformation functor of the $\{\oO_{C_i}(-1)\}$.}
\vspace{0.3cm}

When $Y$ is 3-dimensional but not local, category $\mathscr{A}_{f}$ is the direct sum of categories $\mathscr{A}_{f,C_y}$ where $y$ runs over all closed points in $Y$ with 1-dimensional $f$-fibers $C_y=\bigcup C_{y,i}$. Thus, category $\mathscr{A}_{f}$ \emph{represents} the functor of non-commutative deformations of  $\{\oO_{C_{y,i}}(-1)\}$.

\vspace{0.3cm}
\textbf{Structure of the paper.} In Section \ref{sec_fam_of_sheaves} we study the category $\zZ_A$ of $A$-objects in an arbitrary abelian category $\zZ$. We prove that $\zZ_A$ is equivalent to the category $\Rex(\modu A, \zZ)$ of right exact functors $\modu A \to \zZ$. We say that an $A$-object is \emph{flat} if the corresponding functor is exact.

We show that, for a $\kk$-scheme X and a  $\kk$-algebra $A$, the category $\QCoh(X \times \Spec A)$ is equivalent to $\QCoh(X)_A$ and that sheaves flat over $A$ correspond to flat $A$-objects (see Section \ref{ssec_sheaves_on_prod}). We also discuss a pair of adjoint functors $\zZ_A \leftrightarrow \zZ_B$ induced by a homomorphism $A\to B$.

In Section \ref{ssec_class_def_via_A-obj} we rephrase the classical deformation functor using the language of $A$-objects. This point of view allows us to extend this classical functor to the category $\Bas_1$ (see Section \ref{ssec_cat_Bas_1}). We then consider non-commutative deformations of $n$-objects (see Section \ref{ssec_nc_def_of_n-obj}) as a functor $\Bas_n \to \Sets$.

In Section \ref{sec_abelian_base} we consider abelian categories as bases of  infinitesimal non-commutative flat families. We introduce category $\Df_n$, equivalent to $\Bas_n^{\opp}$, as a subcategory of $\wt{\Df}_n$ (see Section \ref{ssec_Df}). We prove that $\Df_n \subset \wt{\Df}_n$ is the full subcategory of abelian categories with universally bounded Loewy length, see Section \ref{ssec_Df_as_ind-obj}. In the virtue of 
the equivalence $\Df_n \simeq \Bas_n^{\opp}$ we	 define the categorified functor of   non-commutative deformations $\Df_n^{\opp} \to \Sets$ (see Section \ref{ssec_nc_def_Df_to_Set}).

We study subcategories of exact/abelian categories in Section \ref{sec_inf_ext_of_subcat}. We show that extension closure of a full abelian subcategory $\aA \subset \zZ$ is again abelian.

In Section \ref{ssec_def_of_def} we define the deformation functor $\ncDef_\zeta$ of an exact functor $\zeta \colon \aA_0 \to \zZ$.  We discuss ind-objects over $\Df_n$ and, more generally, over $\Fin_{\aA_0}$ in Section \ref{ssec_Fin_as_ind-obj}. We show that if $\zeta$ is fully faithful then $\ncDef_\zeta$ is ind-represented by the extension closure of the full abelian subcategory $\zeta(\aA_0) \subset \zZ$, see Theorem \ref{thm_ind_rep_fully_faith}. In particular, the functor of non-commutative deformations of a \emph{simple collection} $(Z_1,\ldots, Z_n)$ is ind-representable, 
i.e. isomorphic to $\Hom(-, (\cC, Q_{\cC}))$, for some $(\cC, Q_{\cC}) \in \wt{\Df}_n$. 
We also give a detailed example of  non-simple collection for which $\ncDef_{(Z_1,\ldots,Z_n)}$ is not ind-representable (see Section \ref{ssec_obst_to_rep}).

In Section \ref{sec_null_cat_as_nc_def} we study the motivating example related to flopping contractions. In Sections \ref{ssec_A_f_under_base_change} and \ref{ssec_proj_in_A_f} we recall the properties of the null-category $\mathscr{A}_f$. We describe simple objects in $\mathscr{A}_f$ (see Section \ref{ssec_simple_in_Af}) and discuss when $\mathscr{A}_f$ and $\mathscr{A}_{f,C}$ (ind-)represent the functor of non-commutative deformations (see Section \ref{ssec_A_f_C_ind-rep}).

\vspace{0.3cm}
\textbf{Notation.}
We work over an algebraically closed field $\kk$. 

For an algebra $A$, we denote by  $\Modu A$ the category of right $A$-modules. By $\modu A \subset \Modu A$ we denote the full subcategory of finitely presented modules. It is abelian if the algebra $A$ is right coherent. We denote by  $ A\textrm{--Mod}$ the category  of left $A$-modules.

\vspace{0.3cm}
\textbf{Acknowledgements.} 
We are thankful to A. Polishchuk for an useful discussion.
The first named author would like to thank Kavli IPMU and Higher School of Economics for their hospitality. The first named author was partially supported by Polish National Science Centre grants No. 2018/29/B/ST1/01232 and 2021/41/B/ST1/03741. 
The reported study was funded by RFBR and CNRS, project number 21-51-15005.
The work of A. I. Bondal was performed at the Steklov International Mathematical Center and supported by the Ministry of Science and Higher Education of the Russian Federation (agreement no. 075-15-2019-1614).
This research was partially supported by World Premier International Research Center Initiative (WPI Initiative), MEXT, Japan.
This work was supported by JSPS KAKENHI Grant Number JP23K20205.

\section{Families of quasi-coherent sheaves over an affine base}\label{sec_fam_of_sheaves}

We interpret families of quasi-coherent sheaves on a scheme $X$ parameterised by an affine base $\Spec A$, i.e. quasi-coherent sheaves on $X\times \Spec A$, as $A$-objects in $\QCoh(X)$ and as right exact functors $\modu A \to \QCoh(X)$.

\vspace{0.3cm}
\subsection{$A$-objects}~\label{A-objects}\\

Let $\zZ$ be a $\kk$-linear category and $A$ a unital $\kk$-algebra. An $A$-\emph{object} in $\zZ$ is a pair $(Z, \rho)$  of an object $Z\in \zZ$ and a homomorphism $\rho \colon A \to \End(Z)$. We often omit $\rho$.
Denote by $\zZ_A$ the category of $A$-objects in $\zZ$ with morphisms defined in the obvious way.

\begin{PROP}\cite[Proposition B2.2]{ArtZha}
	Category $\zZ_A$ is abelian.
\end{PROP}

For a $\kk$-linear  abelian category $\zZ$ and a homomorphism $\alpha \colon A\to B$,  
we define functor
\begin{align}\label{eqtn_alpha_lower}
&\alpha_*^{\zZ} \colon \zZ_B \to \zZ_A,& &\alpha_*^{\zZ}(Z,\rho_B) = (Z, \rho_B \circ \alpha).&
\end{align}
By formula (\ref{eqtn_alpha_upper}), we shall define its left adjoint $\alpha^*_{\zZ}\colon \zZ_A \to \zZ_B$.

The category of $A$-objects in ${\kk}\textrm{--Mod}$ is $A\textrm{--Mod}$. Then,  $\alpha_*^{{\kk}\textrm{--Mod}}\colon ({\kk}\textrm{--Mod})_B \to ({\kk}\textrm{--Mod})_A$ is just the restriction of scalars functor $\alpha_*\colon B\textrm{--Mod} \to A\textrm{--Mod}$. Similarly, right $A$-modules are $A^{op}$-objects in ${\kk}\textrm{--Mod}$, and the functor assined to $\alpha^{op}:A^{op}\to B^{op}$ is the restriction of scalars $\alpha_*\colon \Modu B \to \Modu A$.

\vspace{0.3cm}
\subsection{$A$-objects and right exact functors out of $\modu A$}\label{ssec_A-objects}~\\

An object $Z$ in a category $\zZ$ yields a functor $\Hom_{\zZ}(Z,-)\colon \zZ\to \textrm{Mod--}\End(Z)$. It has left adjoint, if $\zZ$ is a cocomplete abelian category (i.e. admits arbitrary direct sums):  
\begin{equation}\label{eqtn_general_left_adjoint}
(-) \otimes_{\End{Z}}Z \colon \Modu \End(Z) \to \zZ, 
\end{equation}
defined via presentations of $\End(Z)$-modules by free modules, \emph{cf.} \cite[Theorem 3.7.1]{Pop}.

Let $F\colon \cC\to \dD$ be a functor and $\dD ' \subset \dD$ a full subcategory. We say that $G \colon \dD ' \to \cC$ is a \emph{partial left adjoint} to $F$ if there exist functorial in $C\in \cC$ and $D' \in \dD '$ isomorphisms 
\begin{equation}\label{eqtn_partial_adj}
h_{D',C}\colon \Hom_{\cC}(G(D'),C) \xrightarrow{\simeq}\Hom_{\dD}(D', F(C)).
\end{equation} 
The standard argument shows that the partial left adjoint $G$ is unique up to a unique isomorphism. Moreover, if the categories $\cC$, $\dD'$ and $\dD$ are abelian, then $G$ is right exact.

Now let $A$ be a unital $\kk$-algebra and $(Z, \rho :A\to \End (Z))$ an $A$-object in a $\kk$-linear abelian category $\zZ$. Consider the composite functor 
$$
h_Z(-):= \rho_*\circ\Hom(Z,-)\colon \zZ\to \Modu \End(Z) \to \Modu A
$$
If $(Z, \rho)$ is an $A$-object in an cocomplete abelian $\zZ$, then functor $h_Z$ has left adjoint:
$$
	 (-)\otimes_AZ:=\rho^*(-) \otimes_{\End(Z)}Z\colon \Modu A \to \zZ.
$$ 

If $\zZ$ is not cocomplete, we have a partial adjoint. Recall that an algebra $A$ is right coherent if every finitely generated right ideal in $A$ is finitely presented. If $A$ is right coherent, we denote by $\modu A$ the (abelian) category of finitely presented modules.
\begin{PROP}\label{prop_left_adj_to_Hom}(cf. \cite[Theorem 3.6.3]{Pop})
	\begin{enumerate}
		\item functor $h_Z$
			admits a partial left adjoint $(-)\otimes_A Z\colon \modu A \to \zZ$.
		\item If $A$ is right coherent then $(-)\otimes_AZ$ is characterised as the unique, up to canonical isomorphism, right exact functor $T_Z:\modu A\to \zZ$ such that $T_Z|_{A}=(Z, \rho )$.
	\end{enumerate}  
\end{PROP}
\begin{proof}
	We define functor $T_Z\colon \modu A \to \zZ$, such that $T_Z|_{A} \simeq (Z, \rho )$, and bi-functorial isomorphisms, for any $Z'\in \zZ$ and $M\in \modu A$:
	\begin{align}\label{eqtn_almost_adjun}
		t_{M,Z'} \colon \Hom_\zZ(T_Z(M), Z') \xrightarrow{\simeq}\Hom_A(M, h_Z(Z')).
	\end{align}
	 Then $(-)\otimes_AZ := T_Z(-)$ is the required functor.

	There is a unique additive functor from the category ${\rm Free}(A)$ of finite rank free $A$-modules to $\mathcal{Z}_A$ which sends $A$ to $(Z, \rho)$. Composed with the forgetful functor, it defines $T_Z$ on the full subcategory ${\rm Free} (A) \subset \modu A$.
	For $M\in\modu A$, consider a presentation by free modules of finite rank:
	\begin{equation}\label{eqtn_resolution}
	P_1 \xrightarrow{d} P_0 \to M \to 0
	\end{equation}
	and define $T_Z(M)$ as the cokernel of $T_Z(d)$.
	
	 By \cite[Lemma 6.1]{Pop}, for any morphism $f\colon M \to M'$ and any presentations $P_1 \to P_0 \to M$, $P'_1 \to P'_0 \to M'$, there exist morphisms $f_1$ and $f_0$ such that diagram
	\[
	\xymatrix{P'_1 \ar[r]^{d'} & P'_0 \ar[r] & M' \\
		P_1 \ar[r]^d \ar[u]^{f_1} & P_0 \ar[r] \ar[u]^{f_0} & M \ar[u]^f}
	\]
	commutes. Moreover, the induced morphism $\textrm{coker } T_Z(d) \to \textrm{coker }T_Z(d')$ does not depend on the choice of $f_1$ and $f_0$. In particular, $T_Z(M)$ is defined up to a canonical isomorphism.

	Isomorphism $N\xrightarrow{\simeq} \Hom_A(A,N)$, for any $N\in \modu A$, induces isomorphisms for any $n$:
	\begin{align*}
	&t_{A^{\oplus n}, Z'}\colon \Hom_{\zZ}(Z^{\oplus n}, Z') \to \Hom_A(A^{\oplus n}, \Hom_{\zZ}(Z,Z')).&
	\end{align*}
	
	Then the isomorphism $t_{M,Z'}$ is defined via the commutative diagram
	\[
	\xymatrix{0 \ar[r] & \Hom_A(M, h_Z(Z')) \ar[r] & \Hom_A(P_0, h_Z(Z')) \ar[r] & \Hom_A(P_1, h_Z(Z'))\\
		0 \ar[r] & \Hom_{\zZ}(T_Z(M), Z') \ar[r] \ar[u]^{t_{M,Z'}} & \Hom_{\zZ}(T_Z(P_0), Z') \ar[r] \ar[u]^{t_{P_0,Z'}} & \Hom_{\zZ}(T_Z(P_1), Z')\ar[u]^{t_{P_1,Z'}} }
	\]
	Again, as above, $t_{M,Z'}$ does not depend on the choice of presentation (\ref{eqtn_resolution}).

If $A$ is right coherent then the above construction is a unique extension of $\rho$ to a right exact functor $\modu A \to \zZ$.
\end{proof}
\begin{DEF}\label{def_flat}
We say that $Z\in \zZ_A$ is \emph{flat} if the functor $(-)\otimes_AZ\colon \modu A \to \zZ$ is exact.
\end{DEF}
For abelian categories $\aA$ and $\bB$, we denote by $\Rex(\aA, \bB)$ the category of right exact functors $\aA\to \bB$ and natural transformations thereof.
\begin{PROP}(cf. \cite[Theorem 3.6.2]{Pop})\label{prop_functor_given_by_tensor}
	Let $A$ be a right coherent $\kk$-algebra. The restriction to $A$ and functor $Z\mapsto (-)\otimes_AZ$ define quasi-inverse equivalences: 
	$$ \Rex(\modu A, \zZ) \xlongleftrightarrow{\simeq} \zZ_A.
	$$ 
	These restrict to equivalences between exact functors $\modu A \to \zZ$ and flat $A$-objects.
\end{PROP}
\begin{proof}
	The statement follows from 	Proposition \ref{prop_left_adj_to_Hom}.(2) and Definition \ref{def_flat}. 
\end{proof}
This proposition allows us to think about $A$-objects in a category $\zZ$ as functors $\modu A \to \zZ$, which will be the source of our categorification for the deformation theory.

\begin{EXM}\label{exm_standard_functors_and_A-objects}	
	Given a homomorphism $\alpha \colon A\to B$ of $\kk$-algebras, functor $\alpha_* \colon \modu B \to \Modu A$ is defined by the $B$-object $B_{\alpha}$ in $\Modu A$, i.e. by a $B-A$ bimodule, which is $B$  with the right $A$-module structure induced by $\alpha$: 
\begin{equation}\label{balpha}
	 \alpha_*(-) \simeq (-)\otimes_B B_{\alpha}.
\end{equation}
	 Functor $\alpha_*$ is exact, hence $B_{\alpha}$ is flat as a $B$-object in $\Modu A$.
\end{EXM}

For a cocomplete $\kk$-linear abelian category $\zZ$ and $\alpha \colon A\to B$, we define
\begin{align}\label{eqtn_alpha_upper}
&\alpha^*_{\zZ}\colon \zZ_A \to \zZ_B,& &\alpha^*_{\zZ}(Z, \rho_A) = (B_{\alpha} \otimes_A Z, \rho_B : B \to \End_A(B_{\alpha}) \to \End_{\zZ}(B_{\alpha}\otimes_A Z)),&
\end{align}

The same formula for $\alpha^*_{\zZ}$ is applicable when $\zZ$ is an arbitrary $\kk$-linear abelian category and $B$ is finitely presented as a right $A$-module.

We leave the standard proof of the following fact  to the reader.
\begin{PROP}\label{prop_adjoint_A_obj}
	Consider a $\kk$-linear abelian category $\zZ$  and a homomorphism $\alpha \colon A \to B$ of unital $\kk$-algebras. Assume that $\zZ$ is cocomplete or $B\in \modu A$. Then functor $\alpha^*_{\zZ}$ in (\ref{eqtn_alpha_upper})  is left adjoint to $\alpha_*^{\zZ}$ in  (\ref{eqtn_alpha_lower}).
\end{PROP}

If $Z$ is an $A$-object in $\zZ$,  
functor $\alpha_*(-)\otimes_A Z\colon \modu  B \to \zZ$ is right exact. 
It takes $B\in \modu B$ to $\alpha^*_{\zZ}(Z)$. By Proposition \ref{prop_left_adj_to_Hom}.(2), the unique functor with this property is $(-)\otimes_B \alpha^*_{\zZ}(Z)$. Hence, we have an isomorphism of functors
 \begin{equation}\label{eqtn_alpha_and_tensor}
 (-)\otimes_B \alpha^*_{\zZ}(Z) \simeq \alpha_*(-)\otimes_A Z.
 \end{equation}

\vspace{0.3cm}
\subsection{Fourier-Mukai kernels and $A$-objects in $\QCoh(X)$}\label{ssec_sheaves_on_prod}~\\

We would like to interpret $A$-objects in category $\zZ$ 
as families of objects in $\zZ$ parametrised by a `noncommutative space $\Spec A$'. Here we justify this point of view by geometry, assuming $A$ is commutative and $\zZ$ is the category of (quasi)coherent sheaves on a scheme.

Let $A$ be a Noetherian commutative unital $\kk$-algebra. For a $\kk$-scheme $X$, denote by 
\[
\xymatrix{&X\times \Spec A \ar[dl]_{p_X} \ar[dr]^{p_A}& \\ X & &\Spec A}
\]
the canonical projections.

For $F\in \QCoh(X \times \Spec A)$, the  \emph{abelian Fourier-Mukai} functor is defined as
\begin{equation*} 
\textrm{FM}_X(F)(-) = {p_X}_*(F\otimes p_A^*(-)) \in  \Rex(\modu A, \QCoh(X)).
\end{equation*} 
Here we identify $A$-modules with quasi-coherent sheaves on $\Spec A$.

Using  the identification of $\Rex$ with $A$-objects in Proposition \ref{prop_functor_given_by_tensor}, we define functor 
$$
\mathfrak{F}_X\colon \QCoh(X\times \Spec A) \to \QCoh(X)_A.
$$ 
It takes $F \in \QCoh(X\times \Spec A) $ to the restriction of  $\textrm{FM}_X(F)$ to $A=\oO_{\Spec A}$: 
\begin{equation}\label{eqtn_FM_functor} 
\mathfrak{F}_X(F)= \textrm{FM}_X(F)|_{A}=p_{X*}(F),
\end{equation} 
where $p_{X*}(F)$ has the structure of module over $A_X:=p_{X*}(\oO_{X\times \Spec A})=A\otimes_k \oO_X$.

According to Proposition \ref{prop_functor_given_by_tensor}, functor $(-)\otimes_A\mathfrak{F}_X(F)$ is but 
$\textrm{FM}_X(F)$.
\begin{PROP}\label{prop_sheaves_on_prod_right_ex_fun}
	\begin{enumerate}
		\item  $\mathfrak{F}_X$ is an equivalence of categories,
			\item $\mathfrak{F}_X$ restricts to an equivalence of the full subcategory of sheaves on $X\times \Spec A$ flat over $\Spec A$ and the category of flat $A$-objects in $\QCoh(X)$.
	\end{enumerate}
	\end{PROP} 
\begin{proof} 
	Morphism $p_X$ is affine, hence $p_{X*}$ is an equivalence $\QCoh(X\times \Spec A)\xrightarrow{\simeq} \QCoh(X, p_{X*}\oO_{X\times \Spec A}) = \QCoh(X, A\otimes_k \oO_X)$, \cite[Proposition 1.4.3]{EGAII}. The latter category is equivalent to $\QCoh(X)_A$. Indeed, one checks locally that, for a commutative $k$-algebra $R$, the category ${\rm Mod}(A\otimes_k R)$ is equivalent to the category of $A$-objects in ${\rm Mod}(R)$, \emph{cf} \cite[Example B2.1]{ArtZha}. The equivalence maps an $A\otimes_k R$-module $M$ to itself with the $R$-module structure induced by the homomorphism $R \to A\otimes_k R$, $r\mapsto 1 \otimes r$. The structure of an $A$-object, $\rho \colon A\to \End_R(M)$, is $\rho(a)(m) = (a\otimes 1) m$.

A sheaf $F\in \QCoh(X\times \Spec A)$ is flat over $\Spec A$ if and only if the functor $F\otimes p_A^*(-)$ is exact. As $p_X$ is affine, ${p_X}_*$ is exact and has no kernel, i.e. ${p_X}_*(G)\simeq0$ implies $G\simeq 0$. Therefore, the exactness of $F\otimes p_A^*(-)$ is equivalent to the exactness of $\textrm{FM}_X(F)$. 
\end{proof} 

\begin{PROP}\label{prop_push_forw_of_A-objects}
	Let $X$ be a $\kk$-scheme and $\alpha \colon A \to B$ a homomorphism of Noetherian commutative $\kk$-algebras. Let
	$\alpha_X\colon X \times \Spec B \to X\times \Spec A$ be the morphism induced by $\alpha$. Then the following diagram commutes:
	\begin{equation}\label{eqtn_pushforw_for_objects}
	\xymatrix{\QCoh(X\times \Spec A) \ar[rr]^{\mathfrak{F}_X}&& \QCoh(X)_A\\
	\QCoh(X\times \Spec B) \ar[rr]^{\mathfrak{F}_X} \ar[u]^{{\alpha_X}_*} && \QCoh(X)_B. \ar[u]^{\alpha_*^{\QCoh(X)}}}
	\end{equation}
\end{PROP} 
\begin{proof} 
The statement can be checked locally on $X$. Consider the case $X = \Spec R$. For $M \in {\rm Mod}(B \otimes_k R)$, both $\mathfrak{F}_X \circ \alpha_{X*}(M)$ and $\alpha_*^{\QCoh(X)} \circ \mathfrak{F}_X(M)$ have the $R$-module structure induced by the canonical homomorphism $R\to B \otimes_k R$. Further, $A$-object structures coincide as they are induced by two homomorphisms $A\to B\otimes_k R$ in the commutative diagram:
\[
\xymatrix{A \otimes_k R \ar[d]_{\alpha \otimes \Id} & & A\ar[ll]_{{\rm can}} \ar[d]^\alpha \\
B\otimes_k R && B \ar[ll]_{{\rm can}}} 
\]
\end{proof}

In view of Proposition \ref{prop_push_forw_of_A-objects}, we can interpret $\alpha_*^{\QCoh(X)}$ as a non-commutative version of the push-forward of Fourier-Mukai kernels. Its left adjoint $\alpha^*_{\QCoh(X)}$ acquires the meaning of a non-commutative pull-back of FM kernels.

\begin{COR}\label{cor_pull_back_of_A_obj}
		Let $X$ be a $\kk$-scheme and $\alpha \colon A \to B$ a homomorphism of Noetherian commutative $\kk$-algebras. Let
	$\alpha_X\colon X \times \Spec B \to X\times \Spec A$ be the morphism induced by $\alpha$. Then the following diagram commutes:
	\begin{equation*}
	\xymatrix{\QCoh(X\times \Spec B) \ar[rr]^{\mathfrak{F}_X}&& \QCoh(X)_B\\
		\QCoh(X\times \Spec A) \ar[rr]^{\mathfrak{F}_X} \ar[u]^{{\alpha_X}^*} && \QCoh(X)_A. \ar[u]^{\alpha^*_{\QCoh(X)}}}
	\end{equation*}

\end{COR}
\begin{proof}
	This is implied by Propositions \ref{prop_push_forw_of_A-objects}, \ref{prop_adjoint_A_obj} and the uniqueness of the left adjoint. 
\end{proof}

\section{Non-commutative deformations over the category $\Bas_n$}\label{sec_nc_def_over_Bas_n}

We reformulate the classical deformation theory for sheaves on a $\kk$-scheme $X$ via $A$-objects in the abelian category $\QCoh(X)$. This leads us to flat deformation of sheaves over a non-commutative local Artinian algebra $A$. By extending this viewpoint to deformations of $n$-tuples of objects  in an abelian category $\zZ$ over a basic Artinian algebra $A$ with $n$ isomorphism classes of simple modules, we arrive at a functor $\textrm{ncDef}_{Z_\bcdot}\colon \Bas_n \to \Sets$. 

\vspace{0.3cm}
\subsection{The classical deformation functor via flat $A$-objects}\label{ssec_class_def_via_A-obj}~\\

For a local $\kk$-algebra $A$ we denote by $q_A \colon A\to A/\mathfrak{m}_A \simeq \kk$ the quotient by the maximal ideal. 
By $\textrm{Art}_{\kk}$ we denote the category of finite dimensional  local commutative 	$\kk$-algebras.

Let $F_0$ be a quasicoherent sheaf on a $\kk$-scheme $X$. We define the deformation functor $\textrm{Def}_{F_0} \colon \textrm{Art}_{\kk} \to \Sets$ as follows (see notation in (\ref{eqtn_alpha_upper})):
\begin{equation}\label{eqtn_classical_def_revisited}
\textrm{Def}_{F_0}(A)  = \{(F_A, \f_A)\,|\, F_A\in \QCoh(X)_A \textrm{ a flat  }A \textrm{-object, }  \f_A \colon q^*_{A,\QCoh(X)} F_A \xrightarrow{\simeq} F_0\}/\sim
\end{equation}
where, by definition, $(F_A, \f_A) \sim (F'_A, \f'_A)$ if there exists an isomorphism $\psi \colon F_A\xrightarrow{\simeq} F'_A$ of $A$-objects such that $\f'_A \circ q_{A,\QCoh(X)}^*(\psi) = \f_A$.

For an $\alpha \in \Hom_{\textrm{Art}_k}(A, B)$, the equality $q_B \circ \alpha = q_A$  implies an isomorphism $ q^*_{A,\QCoh(X)} \simeq q^*_{B,\QCoh(X)}\circ \alpha^*_{\QCoh(X)}$.
It allows us to define 
$\textrm{Def}_{F_0}(\alpha) \colon \textrm{Def}_{F_0}(A) \to \Def_{F_0}(B)$ as
\begin{align*}
\Def_{F_0}(\alpha)(F_A, \f_A) = (\alpha^*_{\QCoh(X)}F_A, q^*_{B,\QCoh(X)}\alpha^*_{\QCoh(X)}F_A \xrightarrow{\simeq} q^*_{A,\QCoh(X)} F_A \xrightarrow{\f_A} F_0).
\end{align*}

In view of Proposition \ref{prop_sheaves_on_prod_right_ex_fun}, we recover
the classical deformation functor.
\begin{PROP}\label{prop_2_def_agree}
	Let $F_0$ be a quasicoherent sheaf on $X$. For $A\in \textrm{Art}_{\kk}$, denote by $i_A \colon X\to X\times \Spec A$ the embedding. Then $\Def_{F_0}$ is the classical deformation functor \cite{Schle}:
	\begin{equation}\label{eqtn_classical_def}
	\textrm{Def}_{F_0}(A)  = \{(F, \f \colon i_A^*F \xrightarrow{\simeq} F_0)\,|\, F\in \QCoh(X\times \Spec A) \textrm{ is flat over }A \}/\sim 
	\end{equation}
	where $(F, \f) \sim (F', \f')$ if there exists an isomorphism $\psi \colon F\xrightarrow{\simeq} F'$ such that $\f' \circ i_A^*(\psi) = \f$.
\end{PROP}

\vspace{0.3cm}
\subsection{Category $\Bas_1$ and noncommutative deformations of a single object}\label{ssec_cat_Bas_1}~\\

For a $\kk$-algebra $A$, we denote by $A^\times$ the set of invertible elements in $A$.

Recall that a finite dimensional $\kk$-algebra $A$ is \emph{basic} if $A/\textrm{rad} A$ is isomorphic to a product of copies of $\kk$.
Denote by $\textrm{ncArt}_{\kk}$ the category of finite dimensional basic local $\kk$-algebras. For $A\in {\rm ncArt}_{\kk}$ denote by $q_A \colon A\to A/{\rm rad}A$ the canonical morphism.

We define the category $\Bas_1$. Objects are the same as in $\textrm{ncArt}_{\kk}$, but morphisms are conjugacy classes of homomorphisms:
\begin{equation}\label{eqtn_morph_in_Bas_1}
\begin{aligned} 
&\Hom_{Bas_1}(A, B) = \Hom_{\textrm{ncArt}_{\kk}}(A, B)/\sim , &\\
&\alpha \sim \beta \, \Leftrightarrow \, \exists u\in B^\times \textrm{ such that } q_A(u) =1 \, {\rm and}\, \forall a\in A\, \alpha(a) = u^{-1} \beta(a) u.&
\end{aligned} 
\end{equation}
We say that homomorphisms $\alpha,\beta \colon A\to B$ such that $\alpha = u^{-1}\beta u$ are \emph{conjugate by $u$}. 
One can easily check that the composition of morphisms in $\Bas_1$ is well-defined.

 By analogy with (\ref{eqtn_classical_def_revisited}), for an object $Z_0$ in a $\kk$-linear abelian category $\zZ$,  we define the functor $\ncDef_{Z_0} \colon \Bas_1 \to \Sets$. For $A,B\in \Bas_1$ and $\alpha \in \Hom_{Bas_1}(A, B)$, it is:
\begin{equation}\label{eqtn_nc_def_of_one_obj}
\begin{aligned} 
&\textrm{ncDef}_{Z_0}(A)  = \{(Z, \f)\,|\, Z\in \zZ_A \textrm{ a flat }A \textrm{-object, } \f\colon q^*_{A,\zZ}Z \xrightarrow{\simeq} {Z_0} \}/\sim&\\
&\ncDef_{Z_0}(\alpha)(Z, \f_A) = (\alpha^*_{\zZ}Z, q^*_{B,\zZ}\alpha^*_{\zZ}Z \xrightarrow{\f^{\alpha}}Z_0),&
\end{aligned}
\end{equation}
where $\f^{\alpha}$ is the composite $q^*_{B,\zZ}\alpha^*_{\zZ}Z \xrightarrow{\kappa^\alpha_{Z}} q^*_{A,\zZ} Z\xrightarrow{\f} Z_0$
 and $\kappa^\alpha\colon q_{B, \zZ}^*\alpha^*_{\zZ} \xrightarrow{\simeq} q_{A, \zZ}^*$ is induced by equality $q_B\circ \alpha = q_A$.
By definition, $(Z, \f) \sim (Z', \f')$ if there is an isomorphism $\psi \colon Z\xrightarrow{\simeq} Z'$ of $A$-objects such that $\f' \circ q_{A,\zZ}^*(\psi) = \f$.

The following lemma ensures that the functor $\ncDef_{Z_0}$ is well-defined.

\begin{LEM}\label{lem_conj_morph_iso_families}
	Consider homomorphisms $\alpha, \beta \colon A \to B$ in $\Bas_1$ conjugate by $u\in B^\times$ with $q_B(u) =1$.
	\begin{enumerate}
		\item 
		 $u$ 
		induces an isomorphism of functors $\lambda_u \colon \alpha^*_{\zZ} \xrightarrow{\simeq} \beta^*_{\zZ}$.
		\item
		$\ncDef_{Z_0}(\alpha) = \ncDef_{Z_0}(\beta)$.
	\end{enumerate}
\end{LEM} 
\begin{proof}
	The map $b\mapsto bu^{-1}$ is an isomorphism of $B-A$ bimodules $\psi_u \colon B_\alpha \xrightarrow{\simeq} B_\beta$ (see (\ref{balpha})).
	By (\ref{eqtn_alpha_upper}), it induces an isomorphism of functors $\lambda_u \colon \alpha^*_{\zZ} \xrightarrow{\simeq}\beta^*_{\zZ}$. 
	
	For $(Z, \f) \in \ncDef_{Z_0}(A)$, denote by $\la_{u,Z}\colon \alpha^*_{\zZ}(Z) \xrightarrow{\simeq} \beta^*_{\zZ}(Z)$ the induced isomorphism. To prove part (2), we need to check that $\f_A^{\alpha}=\f_A^{\beta}\circ q^*_{B,Z}(\lambda_{u,Z})$.
	This would follow from 
	the commutativity of the following diagram:
	\begin{equation}\label{eqtn_diag_1}
	\xymatrix{q_{B,\zZ}^*\beta_{\zZ}^*Z \ar[rr]^{\kappa^{\beta}_Z}&& q_{A,\zZ}^*Z \\
	q_{B,\zZ}^*\alpha_{\zZ}^*Z \ar[rru]_{\kappa^{\alpha}_Z}\ar[u]^{q_{B,\zZ}^*(\la_{u,Z})}&& }
	\end{equation}
	Formula (\ref{eqtn_alpha_and_tensor}) for $B=\kk$ implies isomorphisms 
	\begin{align*}
			q_{B,\zZ}^* \beta_{\zZ}^*(Z) \simeq \kk \otimes_{\kk} (\kk_{q_{B}\beta}\otimes_A Z)
	 \xrightarrow{\simeq} (\kk \otimes_{\kk} \kk_{q_{B}\beta})\otimes_A Z \simeq \beta_*{q_B}_*(\kk)\otimes_A Z\\
	 q_{B,\zZ}^* \alpha_{\zZ}^*(Z)\simeq \kk \otimes_{\kk} (\kk_{q_{B}\alpha}\otimes_A Z)
	 \xrightarrow{\simeq} (\kk \otimes_{\kk} \kk_{q_{B}\alpha})\otimes_A Z \simeq \alpha_*{q_B}_*(\kk)\otimes_A Z\\
	 q_{A,\zZ}^*Z \simeq \kk \otimes_{\kk} (\kk_{q_{A}}\otimes_A Z)
	 \xrightarrow{\simeq} (\kk \otimes_{\kk} \kk_{q_{A}})\otimes_A Z \simeq {q_A}_*(\kk)\otimes_A Z.
	\end{align*}
	As morphisms in \eqref{eqtn_diag_1} are induced by morphisms of the `middle' terms $\kk_{q_B\beta}$, $\kk_{q_B \alpha}$ and $\kk_{q_A}$, commutativity of diagram (\ref{eqtn_diag_1})  will follow from the commutativity of the diagram of $A$-modules:
	\begin{equation}\label{eqtn_diag_2}
	\xymatrix{ \beta_*{q_B}_*(\kk) \ar[rr]^{\psi^{\beta}_{\kk}} && {q_A}_*(\kk) \\
	\alpha_*{q_B}_*(\kk) \ar[rru]_{\psi^{\alpha}_{\kk}} \ar[u]^{\tau_{u,{q_B}_*(\kk)}} && }
	\end{equation}
	where,  $\tau_u \colon \alpha_* \xrightarrow{\simeq}\beta_*$ is induced by $\psi_u\colon B_{\alpha}\to B_{\beta}$ (see Example \ref{exm_standard_functors_and_A-objects}) while $\psi^\alpha\colon \alpha_*{q_B}_* \xrightarrow{=} {q_A}_*$ and $\psi^\beta\colon \beta_*{q_B}_* \xrightarrow{=} {q_A}_*$ are identity natural transformations induced by the equalities $q_B\circ \alpha = q_A = q_B\circ \beta$. 
	
By definition $\alpha_*{q_B}_*(\kk) = \kk \otimes_{\kk} \kk_{q_B}\otimes_B B_{\alpha}$ and similarly for $\beta_*{q_B}_*(\kk)$. The morphism $\tau_{u, {q_B}_*(\kk)}$ is then the multiplication with $q_B(u) =1$. As $\psi_{\kk}^\alpha$ and $\psi_{\kk}^\beta$ are identity natural transformations, the commutativity of \eqref{eqtn_diag_2} follows.
	\end{proof}

Note that, for commutative local Artinian algebras $A$, $B$, the map $\Hom_{\textrm{Art}_{\kk}}(A, B) \to \Hom_{\Bas_1}(A, B)$ is an isomorphism, i.e. $\textrm{Art}_{\kk}$ can be considered as a full subcategory of $\Bas_1$. In view of Section \ref{ssec_class_def_via_A-obj}, the restriction of the non-commutative deformation functor to the category of commutative Artinian algebras gives the classical deformation functor:
$$
\ncDef_{F_0}|_{\textrm{Art}_{\kk}} \simeq \Def_{F_0}.
$$
\begin{REM}\label{rem_conj_class}
	Assume that $\Def_{F_0}$ is representable, i.e. $\Hom(A,-)\xrightarrow{\simeq}\ncDef_{Z_0}(-)$, and let $(U,\upsilon)\in \ncDef_{Z_0}(A)$ be the universal family, i.e. the image of $\Id_A$. Then, for any $B\in \Bas_1$ and any element $(Z, \f) \in \ncDef_{Z_0}(B)$, there exists a unique morphism $\alpha \colon A \to B$ in $\Bas_1$ such that $(Z,\f) = \ncDef_{Z_0}(\alpha)(U, \upsilon)$. On the other hand, Lemma \ref{lem_conj_morph_iso_families} shows that, for conjugate homomorphisms $\alpha, \beta \colon A\to B$, we have $\ncDef_{Z_0}(\alpha) (U, \upsilon) = \ncDef_{Z_0}(\beta)(U, \upsilon)$. Thus,
	the equality of conjugate morphisms in $\Bas_1$ is necessary for the representability of the non-commutative deformation functor.
\end{REM}

\vspace{0.3cm}
\subsection{Category $\Bas_n$ and noncommutative deformations of $n$ objects}\label{ssec_nc_def_of_n-obj}~\\

Definition (\ref{eqtn_nc_def_of_one_obj}) of non-commutative deformations of a single object admits a natural generalisation to the deformation of an $n$-tuple $Z_\bcdot=(Z_1,\ldots, Z_n)$ of objects in $\zZ$. 

To this end, we introduce the category $\Bas_n$ of basic Artinian algebras with $n$ maximal ideals. Concretely, objects of $\Bas_n$ are pairs $(A, q_A \colon A \to \kk^{\oplus n})$ of a finite dimensional $\kk$-algebra $A$ and an algebra homomorphism $q_A$ which induces an isomorphism $A/\textrm{rad} A \xrightarrow{\simeq} \kk^{\oplus n}$. Morphisms $(A, q_A) \to (B, q_B)$ in $\Bas_n$ are conjugacy classes of $\kk$-algebra homomorphisms $\alpha \colon A \to B$ such that $q_B \circ \alpha = q_A$. As before, $\alpha, \beta \colon A \to B$ are conjugated if there exists $u\in B^\times$, $q_B(u)=1$ such that $\alpha = u^{-1}\beta u$.

Given $n$-tuple $Z_\bcdot$,
we consider  $\bigoplus_{i=1}^nZ_i$ as a $\kk^{\oplus n}$-object.
Define the non-commutative deformation functor $\ncDef_{Z_\bcdot} \colon \Bas_n \to \Sets$ by:
\begin{equation}\label{eqtn_nc_def_n_obj}
\begin{aligned} 
&\textrm{ncDef}_{Z_\bcdot}(A, q_A)  = \{(Z, \f)\,|\, Z\in \zZ_A \textrm{ a flat }A \textrm{-object, }\f \colon q_{A,\zZ}^*Z \xrightarrow{\simeq} \bigoplus_i Z_i \}/\sim& \\
&\ncDef_{Z_\bcdot}(\alpha)(Z, \f) = (\alpha^*_{\zZ}Z, q_{B,\zZ}^*\alpha^*_{\zZ}Z \xrightarrow{\f^\alpha} \bigoplus Z_i),&
\end{aligned}
\end{equation}
where $\alpha \colon A \to B$ is a morphism in $\Bas_n$ and $\f^\alpha$ is the composite 
$q^*_{B,\zZ}\alpha^*_{\zZ}Z \xrightarrow{\kappa^\alpha_{Z}} q^*_{A,\zZ} Z \xrightarrow{\f} \bigoplus Z_i$ and $\kappa^\alpha\colon q_{B, \zZ}^*\alpha^*_{\zZ} \xrightarrow{\simeq} q_{A, \zZ}^*$ is induced by equality $q_B\circ \alpha = q_A$.
By definition, $(Z, \f) \sim (Z', \f')$ if there exists an isomorphism $\psi \colon Z\xrightarrow{\simeq} Z'$ of $A$-objects such that $\f'\circ q_{A,\zZ}^*( \psi) = \f$.

An argument analogous to the one  Lemma \ref{lem_conj_morph_iso_families}.(2) shows that functor $\ncDef_{Z_\bcdot}$ is well-defined.

\begin{REM} 
	Consider an $n$-tuple $V_\bcdot = (V_1,\ldots,V_n)$ of right modules over a $\kk$-algebra $R$. For $A\in \Bas_n$, an $A$-object in $\Modu R$ is an $A-R$ bimodule. Hence,  an element of $\ncDef_{V_\bcdot}(A)$ is an $A-R$ bimodule $W$ and an isomorphism $\f\colon q_A^*(W) = (A/\textrm{rad}A) \otimes_A W \xrightarrow{\simeq} \bigoplus V_i$. 	Moreover, $W$ is flat as an $A$-module, hence projective (see \cite[Theorem P]{Bass}). 
	
	A choice of splitting of $q_A\colon A \to \kk^{\oplus n}$ corresponds to a choice of orthogonal idempotents $e_1,\ldots, e_n$ in $A$ such that $\sum e_i =1$. Then, $W$ being projective implies that $W\simeq \bigoplus Ae_i \otimes_k W_i$ as a left $A$-module. The isomorphism $\f$ is then induced by  isomorphisms $\f_i \colon\kk \otimes_{\kk} W_i \xrightarrow{\simeq} V_i$. 
	The structure of an $R$-module is given by a homomorphism
	\begin{equation}\label{eqtn_Laudal}
	\eta \colon R\to \End_A(\bigoplus_{i=1}^n Ae_i \otimes_k W_i) =(\Hom_{\kk}(W_i, W_j) \otimes_k e_i A e_j),
	\end{equation}
	where we view the right hand side as a matrix algebra. Finally, isomorphism of $A$-objects is given by an isomorphism of $\bigoplus A e_i \otimes_{\kk} W_i$ which respects the $R$-module structure, i.e. commutes with \eqref{eqtn_Laudal}.

	In \cite{Lau} Laudal defined the deformation functor  $\ncDef^{\textrm{Lau}}$ of an $n$-tuple $V_\bcdot$ of $R$-modules. The functor is defined on the category $\Bas_n^\textrm{Lau}$ of algebras $A$ in $\Bas_n$ with a fixed splitting $\iota$ of $q_A $. The choice of $\iota$ allows to view $A$ as an $n\times n$ matrix $A =(e_i A e_j)$. Laudal considers the left $A$-module $(e_iAe_j\otimes_{\kk} V_j)$ and notes that the choice of $q_A \colon A\to \kk^n$ yields the homomorphism 
	\begin{equation}\label{eqtn_Laudal2}
\wt{q} \colon (\Hom_{\kk}(V_i, V_j) \otimes_{\kk}e_i A e_j) \to \bigoplus_{i} \End_{\kk}(V_i).	
\end{equation}
of endomorphism rings. Moreover, the right $R$-module structure on $V_i$'s is defined by a homomorphism $\eta_0 \colon R\to \bigoplus_{i=1}^n \End_{\kk}(V_i)$. By definition, 
\begin{align*}
	\ncDef^{\textrm{Lau}}(A) = \{\eta \textrm{ as in \eqref{eqtn_Laudal}}\,|\,\wt{q} \circ \eta = \eta_0\}/\sim, 
\end{align*}
where $\eta \sim \eta'$ if they differ by an inner automorphism in $ (\Hom_{\kk}(V_i, V_j) \otimes_{\kk} e_i A e_j)$.

Consider functor $F\colon \Bas_n^\textrm{Lau} \to \Bas_n$ defined by forgetting the splitting of $q_A$ and identifying conjugate homomorphisms of algebras. Our definition of non-commutative deformations agrees with that of Laudal in the sense that the diagram commutes:
	\[
	\xymatrix{ \Bas_n^\textrm{Lau} \ar[rr]^{\ncDef^\textrm{Lau}} \ar[d]_F& & \Sets \\ \Bas_n \ar[urr]_{\ncDef}&&}
	\]
Indeed, let $(A, q_A, \iota)$ be an object of $\Bas_n^\textrm{Lau}$. Consider the map $\ncDef_{V_{\bcdot}}(F(A)) \to \ncDef^{\rm Lau}_{V_{\bcdot}}(A)$ which maps $(W, \f)$ to $\eta$ as in \eqref{eqtn_Laudal}, where $W_i$ are identified with $V_i$ via $\f_i$. Since $A\to A/{\rm rad}A$ is a homomorphism of $A$-modules, after tensoring with $(V_i)$ it becomes a homomorphism of $R$-modules, i.e. $\wt{q} \circ \eta = \eta_0 $ as required. If $\psi \colon W\to W'$ is an isomorphism of $A-R$ bimodules, the corresponding $\eta$ and $\eta'$ differ by conjugation with $\psi$, i.e. the map is well-defined. It is immediate to see that its inverse maps $\eta$ to the $A$-module $\bigoplus Ae_i \otimes_{\kk} V_i$ with the $R$-module structure given by $\eta$ and isomorphism $\f $ induced by $A/{\rm rad} A\simeq \kk^{\oplus n}$.
	
		Since functor $F$ is not an equivalence, as it identifies conjugate morphisms,  Laudal's deformation functor differs from the one considered in this paper, in particular, it is almost never ind-representable (\textit{cf.} section \ref{ssec_ind_rep_ff}).
	
	\emph{Warning:} In the notation of Laudal \cite{Lau} the roles of $A$ and $R$ are exchanged.
\end{REM}

\section{An abelian category as the base of a non-commutative deformation}\label{sec_abelian_base}

We introduce the category $\wt{\Df}_n$ of Hom and $\Ext^1$-finite $\kk$-linear abelian categories with $n$ isomorphism classes of simple objects and its full subcategory $\Df_n$ of Deligne finite categories. 
We show that $A \mapsto \modu A$ yields an equivalence $\Bas_n \simeq \Df_n^{\opp}$. 
As the central point of categorification in this paper, we reinterpret  the functor of non-commutative deformations as a functor $\Df_n^{\opp} \to \Sets$.

\vspace{0.3cm}
\subsection{An anti-equivalence between $\Bas_n$ and $\Df_n$}\label{ssec_Df}~\\

Recall that an object of an abelian category is \emph{simple} if it has no proper subobjects. An object is \emph{semi-simple} if it is isomorphic to a finite direct sum of simple objects. For an abelian category $\aA$, we denote by $\sS(\aA)$ the full subcategory of semi-simple objects.

We say that a $\kk$-linear abelian category is \emph{of finite length} if it is $\Hom$ and $\Ext^1$-finite, contains finitely many non-isomorphic simple objects, and any object admits a finite filtration with semi-simple graded factors. An abelian category of finite length is a \emph{Deligne finite category} if it has a projective generator.
By \cite[Proposition 2.14]{Del}, a Deligne finite category $\aA$ with projective generator $P$ is equivalent to the category of finitely generated modules over the (finite dimensional) $\kk$-algebra $\End_{\aA}(P)$. 

We define the category $\wt{\Df}_n$. Its objects are pairs $(\aA, Q_\aA)$, where $\aA$ is a $\kk$-linear abelian category of finite length and $Q_\aA \colon \modu \kk^{\oplus n} \to \aA$ is a fully faithful functor with image $\sS(\aA)$. For $(\aA, Q_{\aA})$, $(\bB, Q_{\bB})\in \wt{\Df}_n$,
\begin{align*}
\Hom_{\wt{\Df}_n}((\aA, Q_{\aA}), (\bB, Q_{\bB}) ) = \{(F, \f \colon F\circ Q_{\aA} \xrightarrow{\simeq} Q_B)\,|\, F \colon \aA \to \bB \textrm{ is exact }\}/\sim
\end{align*}
where $(F,\f) \sim (F', \f')$ if there exists $\psi \colon F\to F'$ such that $\f' (\psi\circ Q_{\aA}) = \f$.
One easily checks that the composition in $\wt{\Df}_n$ induced by composition of exact functors is well-defined.
We often omit $Q_{\aA}$ and write $\aA\in \wt{\Df}_n$.

By $\Df_n$ we denote the full subcategory of $\wt{\Df}_n$ consisting of Deligne finite categories.

Consider functor $\Upsilon \colon \textrm{Bas}_n^{\opp} \to \Df_n$ defined on $A\in \textrm{Bas}_n$ and $\alpha \in \Hom_{\textrm{Bas}_n}(A,B)$ as follows:
\begin{align}\label{eqtn_Upsilon}
 &\Upsilon (A, q_A) =(\modu A, {q_A}_*),& &\Upsilon(\alpha) = (\alpha_* \colon\modu B \to \modu A,\textrm{can}_\alpha\colon \alpha_* {q_B}_* \xrightarrow{\simeq}{q_A}_*),&
\end{align}
where $\textrm{can}_\alpha$ is the identity natural transformation induced by the equality $q_B \circ \alpha = q_A$.

\begin{LEM}
	Functor $\Upsilon$ is well-defined.
\end{LEM}
\begin{proof}
Since $\textrm{rad} A$ annihilates any semi-simple $A$-module,
	functor ${q_A}_*$ is an equivalence of $\modu \kk^{\oplus n}$ with $\sS (\modu A)$, i.e. $(\modu A, {q_A}_*)\in \Df_n$.
	For $\alpha \in \Hom_{\Bas_n}((A,q_A),(B,q_B))$, functor $\alpha_*$ is exact and $\textrm{can}_\alpha$ is an isomorphism.  
	For homomorphisms $\alpha,\beta\colon A \to B$ conjugated by $u\in B^\times $, $q_B(u) =1$, the map $b\mapsto bu^{-1}$ is an isomorphism of $B-A$ bimodules $\psi_u \colon B_\alpha \xrightarrow{\simeq} B_\beta$. Isomorphism (\ref{balpha}) implies an isomorphism $\psi_u \colon \alpha_*\xrightarrow{\simeq} \beta_*$.Since $q_B(u)=1$, isomorphism ${\psi_u}{{q_B}_*}$ is the identity on the $\kk$-vector spaces underlying $B$-modules in the essential image of ${q_B}_*$. Since the same is true about $\textrm{can}_\alpha$ and $\textrm{can}_\beta$, we have  $\textrm{can}_\alpha \circ {\psi_u}^{-1}{{q_B}_*} = \textrm{can}_\beta$.
\end{proof}

\begin{THM}\label{thm_equiv_of_cat}
	Functor $\Upsilon$ defines an anti-equivalence of categories $\Bas_n$ and $\Df_n$.
\end{THM}
\begin{proof}

First, check that $\Upsilon$ is faithful. Assume $\alpha, \beta \in \textrm{Hom}_{\textrm{Bas}_n}(A,B)$ are such that there exists an isomorphism  $\psi \colon \alpha_*\xrightarrow{\simeq} \beta_*$ with  $\textrm{can}_\beta \circ \psi{{q_B}_*} = \textrm{can}_\alpha$. Then there is
	 	an isomorphism $\lambda\colon B_\alpha \xrightarrow{\simeq }B_\beta$ of $B-A$ bimodules
	 	(see Example \ref{exm_standard_functors_and_A-objects}).
	Any automorphism of the left $B$-module ${}_B B$ is of the form $ b \mapsto bu$, for some  $u \in B^\times$. Since $\lambda$ is an isomorphism of right $A$ modules, we have $\alpha(a) u = u \beta(a)$, for any $a \in A$, i.e. $\alpha = u \beta u^{-1}$.  Functors $\textrm{can}_\alpha$ and $\textrm{can}_\beta$ are identities on the underlying vector spaces, while $\psi{{q_B}_*}$ is the right multiplication with $q_B(u)$. As $\textrm{can}_\beta \circ \psi{{q_B}_*} = \textrm{can}_\alpha$, we get that $q_B(u)=1$.
	
	Next, check that $\Upsilon$ is full. Let $(\Psi \colon  \modu B\to \modu A, \f \colon \Psi \circ q_{B*} \xrightarrow{\simeq} q_{A*})$ be a morphism in $\Df_n$. 
	Functor $\Psi$ is isomorphic to the functor $(-)\otimes_B M$, for the $B- A$ bimodule $M = \Psi(B)$. Indeed, both functors are right exact and their restrictions to $B\subset \modu B$ are isomorphic, hence, by Proposition \ref{prop_left_adj_to_Hom}.(2), $\Psi \simeq (-)\otimes_BM$.
	As $\Psi$ is exact, $M$ is flat as a left $B$ module. Algebra $B$ is Artinian. Then $M$ is projective by \cite[Theorem P]{Bass},  i.e. ${}_B M \simeq \bigoplus_{j=1}^n (P_j^B)^{\oplus b_j}$, where $P_i^B$ is an irreducible projective cover of the simple $B$-module $S_i^B$, the image of the $i$'th simple $\kk^{\oplus n}$-module under the functor $q_{B*}$. Object 
$\Psi(S_i^B) =S_i^B \otimes_B (\bigoplus_j (P^B_j)^{\oplus b_j})$ is isomorphic to  $(S_i^A)^{\oplus b_i}$, i.e. to $b_i$ copies of the $i$'th simple $A$-module. As $\f$ is an isomorphism of $\Psi \circ q_{B*}$ and $q_{A*}$, $\Psi(S_i^B) \simeq S_i^A$. Hence, 
	$b_i =1$, for any $i$.
	 Since $B$ is basic, we conclude that ${}_B M \simeq \bigoplus_i P_i^B \simeq {}_B B$.
	 
	 Isomorphism $\f$ can be viewed as an isomorphism $\kk^{\oplus n}_{q_B} \otimes_B M_A \xrightarrow{\simeq} \kk^{\oplus n}_{q_A}$. We choose $m \in M$ such that $\f(\mathbbm{1}_{\kk^{\oplus n}}\otimes m) = \mathbbm{1}_{\kk^{\oplus n}}$. Under the identification of $M$ with $B$ as a left $B$-module, $m$ is identified with an element in $B$ which is $\mathbbm{1}_B$ modulo the radical, i.e. an invertible element in $B$. Hence, $m$ is a generator of $M$ as a $B$-module.
	 The choice of $m$ gives $\alpha^{\opp} \colon A^{\opp} \rightarrow \textrm{End}_B({}_B M) \simeq \textrm{End}_B({}_B B) = B^{\opp}$. Then $(\Psi, \f)\sim (\alpha_*, {\rm can}_{\alpha})$. Indeed, let $\psi\colon \alpha_* \to \Psi$ be induced by the map $B \to M$ which maps the unit $\mathbbm{1}_B$ to $m$. Then the natural transformation $\f \circ \psi q_{A*}\colon \alpha_*q_{B*} \to q_{A*}$ is induced by the composition of morphisms of $\kk^{\oplus n} -B$ modules
	 \begin{align*}
	 	&\kk^{\oplus n}_{q_B} \otimes_B B_{\alpha} \xrightarrow{\psi q_{A*}} \kk^{\oplus n}_{q_B} \otimes_B M_A \xrightarrow{\f} \kk^{\oplus n}_{q_A},& &\mathbbm{1}_{\kk^{\oplus n}} \otimes \mathbbm{1}_B \mapsto \mathbbm{1}_{\kk^{\oplus n}} \otimes m \mapsto \mathbbm{1}_{\kk^{\oplus n}},
	 \end{align*}
	 i.e. it equals the canonical isomorphism ${\rm can}_\alpha$.

	Finally, check that $\Upsilon$ is essentially surjective. 
	Any category $(\aA, Q_\aA)\in \Df_n$ is equivalent to the category of modules over the endomorphism algebra of a projective generator (\emph{cf.} \cite{Del}). Since any finite dimensional algebra is Morita equivalent to a basic one (c.f \cite[Corollary I.6.3]{ASS}), it follows that $\aA \simeq \modu A$, for some finite dimensional basic $\kk$-algebra $A$. Since $\Upsilon$ is full, exact $Q_\aA\colon \modu \kk^{\oplus n} \to \modu A$ is isomorphic to ${q_A}_*$, for a homomorphism $q_A \colon A \to \kk^{\oplus n}$. Homomorphism $q_A$ is unique, because all invertible elements in $\kk^{\oplus n}$ are central. Since the essential image of $Q_\aA$ is the subcategory of semi-simple objects in $\modu A$, $q_A$ induces an isomorphism $A/\textrm{rad}A \simeq \kk^{\oplus n}$, i.e. $(A, q_A) \in \Bas_n$.
\end{proof}

\vspace{0.3cm}
\subsection{Projective generators for categories with bounded Loewy length}\label{ssec_Df_as_ind-obj}~\\

Let $(\aA, Q_{\aA}) \in \wt{\Df}_n$. Since category $\mathcal{A}$ is both artinian and noetherian, functor $Q_{\aA}$ admits left and right adjoint functors $Q_{\aA}^*, Q_{\aA}^! \colon \aA \to \modu \kk^{\oplus n}$. For $A\in \aA$, $Q_{\aA}Q_{\aA}^*(A)$ is the maximal semi-simple quotient of $A$, while $Q_{\aA}Q_{\aA}^!(A)$ is its maximal semi-simple subobject. We define the \emph{radical} endo-functor $\rad\colon \aA\to \aA$ as the kernel of the $Q_{\aA}^* \dashv Q_{\aA}$ adjunction unit $\varepsilon \colon \Id_{\aA} \to Q_{\aA} Q_{\aA}^*$.
For any $A\in \aA$, the natural transformation $\iota \colon \rad \to \Id_{\aA}$ fits into a short exact sequence
\begin{equation}\label{eqtn_radical_as_kernel1} 
0 \to \rad A \xrightarrow{\iota_A} A \xrightarrow{\varepsilon_A} Q_{\aA}Q_{\aA}^* A \to 0.
\end{equation} 
For $f\colon A \to A'$, $\rad f\colon \rad A \to \rad A'$ is just the restriction of $f$ to the subobject $\rad A$.

Recall that $B \subset A$ is \emph{superfluous} if $B + C = A$ implies $C=A$, for any subobject $C\subset A$.
By \cite[Theorem 2.3]{Koh}, $\rad A\subset A$ is the maximal superfluous subobject in $A$. 

\begin{LEM}\label{lem_R_surj_is_surj1}
	Functor $\rad$ preserves mono- and epimorphisms.
\end{LEM}
\begin{proof}
	Let $\f \colon A_1 \to A_2$ be a monomorphism. As $\iota_{A_2} \circ \rad \f  = \f \circ \iota_{A_1}\colon \rad A_1 \to A_2$ is a monomorphism, so is $\rad \f$.
	
	Now let $\f \colon A_1 \to A_2$ be an epimorphism and K its kernel. Functor $Q_{\aA}$ is exact and $Q^*_{\aA}$ is right exact, hence sequence  
	$$
	Q_{\aA} Q_{\aA}^* K \to Q_{\aA} Q_{\aA}^* A_1 \xrightarrow{Q_{\aA}Q_{\aA}^*(\f)} Q_{\aA} Q_{\aA}^*A_2 \to 0
	$$
	is exact. Let $K'$ be the kernel of $Q_{\aA}Q_{\aA}^*(\f)$, i.e. the image of $Q_{\aA}Q_{\aA}^*K$ in $Q_{\aA}Q_{\aA}^*A_1$. The natural transformation $\varepsilon \colon \Id \to Q_{\aA}Q_{\aA}^*$ applied to the canonical morphism $K \to A_1$ gives the left square of the commutative diagram 
	 \begin{equation}\label{eqtn_dia}
	 	\xymatrix{0 \ar[r] & K' \ar[r]& Q_{\aA} Q_{\aA}^* A_1 \ar[r]^{Q_{\aA} Q_{\aA}^*(\f)} & Q_{\aA} Q_{\aA}^* A_2 \ar[r] & 0 \\
		0 \ar[r] & K \ar[r] \ar[u] & A_1 \ar[r]^{\f} \ar[u] & A_2 \ar[r] \ar[u] & 0} 
	 \end{equation}
	with exact rows.
	
	The composite $K\xrightarrow{\varepsilon_K} Q_{\aA}Q_{\aA}^*K \to K'$ of two surjective morphisms is surjective, hence the snake-lemma applied to \eqref{eqtn_dia} 	shows, in view of sequence (\ref{eqtn_radical_as_kernel1}), that $\rad A_1 \to \rad A_2$ is surjective. 
\end{proof}

Denote by $\aA_l$ the full subcategory of $\aA$ consisting of $A \in \aA$ such that $\rad^{l+1}A= 0$. 

\begin{PROP}\label{lem_R_surjective1}
	The subcategory $\aA_l \subset \aA$ is closed under quotients and subobjects (but, in general, not under extensions). In particular, $\aA_l$ is abelian. 
\end{PROP}
\begin{proof}
	The statement follows from Lemma \ref{lem_R_surj_is_surj1}.
\end{proof}

Recall that a \emph{Loewy length} of an object $A\in \aA$ is the  minimal length $l$ of a filtration $0=F_{l+1} \subset F_l \subset \ldots \subset F_0 =A$ with semi-simple graded factors $F_i/F_{i+1}$.

\begin{LEM}\label{lem_radical_and_filtration}
$\aA_l$ is the subcategory of objects in $\aA$ with Loewy length bounded by $l$.
\end{LEM}
\begin{proof}
	For $A\in \aA$ with $\rad^{l+1}A=0$, the filtration $F_i:=\rad^iA$ is the required one.
	
Now let $	0 = F_{l+1} \subset F_l \subset \ldots \subset F_0 =A$ be a filtration with semi-simple graded factors. Since 
the maximal semi-simple quotient of $A$ is $A/\rad A$, we have $\rad A \subset F_1$. 
By iterating, we get $\rad ^i A\subset F_i$, in particular, $\rad^{l+1}A=0$.
\end{proof}

The following theorem can be deduced from Gabriel's description of abelian categories of finite length \cite{Gabriel}. We give a direct proof based on radicals.

\begin{THM}\label{thm_Df_as_bound_on_Loewy}
	Let $\aA$ be an object in $\wt{\Df}_n$. Then $\aA\in \Df_n$ ($\aA$ admits a projective generator) if and only if $\aA \simeq \aA_l$ for some $l\in \mathbb{N}$ (Loewy length of objects is uniformly bounded).
\end{THM}
\begin{proof}
	Assume that  $\aA=\aA_l$ and check that it admits a projective generator. We proceed by induction on $l$. The case $l=0$ is clear.
	
	Let $P $ be a projective generator for $\aA_{l-1}$ and $(S_1,\ldots, S_n)$ the set of irreducible objects in $\aA$.   
	Consider $\wt{P}$, the universal extension of $P$ by $T:=\bigoplus S_i \otimes \Ext^1(P, S_i)^\vee$:
	\begin{equation}\label{eqtn_univ_ext}
	0 \to T  \xrightarrow{i} \wt{P} \xrightarrow{\pi} P \to 0.
	\end{equation}
	The extension is universal in the sense that for any semisimple $L$ and for any  $\eta \in \Ext^1(P, L)$  there is a unique $\gamma :T\to L$ such that $\eta$ is induced from (\ref{eqtn_univ_ext}) by $\gamma$. In other words, applying functors $\Ext^{\bullet}(-,L)$ to (\ref{eqtn_univ_ext}) for such $L$ implies isomorphism $\Hom (T,L)\simeq\Ext^1(P,L)$.
	 
	 We shall show that $\wt{P}$ is a projective generator for $\aA$.
	First, $\rad^l\wt{P}=T$. Indeed, morphism $\pi |_{\rad^l\wt{P}}=\rad^l\pi\colon \rad^l\wt{P} \to \rad^lP =0$ is zero, hence  $\rad^l\wt{P} \subset T$. If $\rad^l\wt{P}\ne T$, then the extension  
	\begin{equation}\label{eqtn_ext}
		0 \to T/ \rad^l\wt{P}\to \wt{P}/\rad^l\wt{P}\to P \to 0
	\end{equation}
 induced from (\ref{eqtn_univ_ext}) by the morphism $T\to T/ \rad^l\wt{P}$ is non-trivial by the universal property. 
	 Clearly, $\rad^l(\wt{P}/\rad^l\wt{P} )=0$, hence
	 $\wt{P}/\rad^l\wt{P} \in \aA_{l-1}$. Also, as $T$ is semi-simple, it is an object of $\aA_0\subset \aA_{l-1}$, hence so is its quotient $T/R^l\wt{P}$. Then extension \eqref{eqtn_ext}  contradicts the fact that $P$ is projective in $\aA_{l-1}$. It implies that  $ \rad^l\wt{P}=T$.
	 
	 	Next, we check that $\Ext^1_{\aA}(\wt{P}, S) =0$, for any simple $S\in \aA$. An extension 
	 	\begin{equation}\label{eqtn_extension}
	 	0 \to S \xrightarrow{\alpha} M\xrightarrow{\beta} \wt{P} \to 0
	 	\end{equation} 
	 	yields 
	 	a commutative diagram with exact rows and columns:
	 	\[
	 	\xymatrix{0 \ar[r]& P \ar[r]^{\Id} & P \\ 
	 		S \ar[r]^{\alpha} \ar[u] & M \ar[r]^{\beta} \ar[u]^g & \wt{P} \ar[u]^{\pi} \\
	 		S \ar[r] \ar[u]^{\Id} & K \ar[r] \ar[u] & T \ar[u]^i}
	 	\]
In the long exact sequence obtained by applying $\Ext^{\bullet}(-,S)$ to (\ref{eqtn_univ_ext}) :
	 	$$
	 	\Hom(T ,S) \to \Ext^1(P, S) \to \Ext^1(\wt{P}, S) \xrightarrow{\tau} \Ext^1(T,S)
	 	$$
	 	the first map is an isomorphism by the universal property. Hence $\tau$ is a monomorphism. In other words, if (\ref{eqtn_extension}) was non-trivial, so was the extension of $T$ by $S$. If this was the case, then $T$ is the maximal semi-simple quotient of $K$. Indeed, the only quotient of $K$ with an epimorphism to $T$ is $K$ itself which is not semi-simple as the extension of $T$ by $S$ is non-trivial. It follows that  $S = \rad K$ is a superfluous subobject in $K$.

Since $\rad^{l}P =0$, the object $\rad^lM$ is contained in the kernel $K$ of $g$.  Lemma \ref{lem_R_surj_is_surj1} implies that $\rad^l\beta \colon \rad^lM \to \rad^l\wt{P}= T$ is surjective. Then $S + \rad^lM\subset K$ contains $S$ and admits a surjective map to $T$, hence $S + \rad^lM =K$. As $S\subset K$ is superfluous, $\rad^lM =K$. Then $\rad^{l+1}M = \rad K =S$ is non-zero, which contradicts the assumption $\aA = \aA_l$. Hence, (\ref{eqtn_extension}) splits and $\Ext^1_{\aA}(\wt{P},S) =0$. As $S\in \aA$ was an arbitrary simple, $\wt{P} \in \aA$ is projective. 

As $P$ has a non-zero morphism to any $S_i$, so does $\wt{P}$, hence $\wt{P}$ is a  generator for $\aA$.

Conversely, let $P\in \aA$ be a projective generator  and $R^{l+1}(P)  =0$, i.e. $P \in \aA_l$.  Every $M \in \aA$ is a quotient of a direct sum of copies of $P$. As the subcategory $\aA_l \subset \aA$ is closed under direct sums and quotients by Proposition \ref{lem_R_surjective1}, it follows that $M$ is an object of $\aA_l$.
\end{proof}

\vspace{0.3cm}
\subsection{The categorified non-commutative deformation functor}\label{ssec_nc_def_Df_to_Set}~\\

Let $Z_\bcdot=(Z_1,\ldots,Z_n)$ be an $n$-tuple of objects in a $\kk$-linear abelian category $\zZ$. The structure of $\kk^{\oplus n}$-object on $\bigoplus_i Z_i$ 
yields a functor (see Proposition \ref{prop_functor_given_by_tensor}):
\begin{equation}\label{eqtn_functor_zeta}
\zeta \colon \modu \kk^{\oplus n} \to \zZ.
\end{equation}

We define a functor $\ncDef_{\zeta} \colon \Df_n^{\opp} \to \Sets$ via:
\begin{equation}\label{eqtn_ncDef_cat}
\ncDef_{\zeta}(\aA, Q_{\aA}) = \{(T,\tau\colon T\circ Q_{\aA}\xrightarrow{\simeq} \zeta )\,|\, T\colon \aA\to \zZ \textrm{ is exact}\}/\sim
\end{equation}
where $(T,\tau) \sim (T',\tau')$ if there exists $\psi \colon T\xrightarrow{\simeq}T'$ such that $\tau'\circ \psi_{ Q_\aA} = \tau$. Note that exactness of the functor $T$ in \eqref{eqtn_ncDef_cat} corresponds to flatness of the $A$-object  in \eqref{eqtn_nc_def_n_obj}, see Proposition  \ref{prop_functor_given_by_tensor}.

For $(F,\f)\in \Hom_{\Df_n}(\aA, \bB)$, we 
define $\ncDef_\zeta(F,\f) \colon \ncDef_{\zeta}(\bB) \to \ncDef_{\zeta}(\aA)$ via:
\begin{equation}
\ncDef_\zeta(F,\f)(T,\tau) = (T\circ F, T\circ F \circ Q_\aA \xrightarrow{T\f} T\circ Q_\bB \xrightarrow{\tau} \zeta).
\end{equation}

\begin{PROP}\label{prop_two_def_coincide}
For an $n$-tuple $Z_\bcdot=(Z_1,\ldots,Z_n)$ of objects of a $\kk$-linear abelian category $\zZ$ and $\zeta = T_{\bigoplus_i Z_i} \colon \modu \kk^{\oplus n} \to \zZ$ the composite $\Bas_n \xrightarrow{\Upsilon} \Df_n^{\opp} \xrightarrow{\ncDef_\zeta} \Sets$ is isomorphic to $\ncDef_{Z_\bcdot}$.
\end{PROP}
\begin{proof}
	Let $(Z_A, \f_A\colon q_A^* Z_A \to\bigoplus_i Z_i)$ be an element of $\ncDef_{Z_\bcdot}(A,q_A)$. By Proposition \ref{prop_functor_given_by_tensor}, the flat object $Z_A$ uniquely determines an exact functor $(-)\otimes_A Z_A\colon \modu A \to \zZ$. Isomorphism (\ref{eqtn_alpha_and_tensor}) implies that $\f_A$ induces an isomorphism ${q_A}_*(-)\otimes_A Z_A \simeq (-) \otimes_{\kk^{\oplus n} } q_{A,\zZ}^*(Z_A) \xrightarrow{\simeq} (-)\otimes_{\kk^{\oplus n}} \bigoplus_i Z_i \simeq \zeta$. An
	isomorphism $\theta \colon Z_A \xrightarrow{\simeq}Z'_A$ yields a unique isomorphism $(-)\otimes_A Z_A \xrightarrow{\simeq}\otimes_A Z'_A$ (see Proposition \ref{prop_functor_given_by_tensor}). 
	
	One easily checks that the constructed isomorphism $\ncDef_{Z_\bcdot}(A, q_A) \to \ncDef_{\zeta} \circ \Upsilon(A,q_A)$ extends to a natural isomorphism of functors.
\end{proof}

\section{Infinitesimal extensions of subcategories}\label{sec_inf_ext_of_subcat}

The Gabriel product of two strictly full subcategories of an exact category allows us to define the $k$-infinitesimal extension $\ngbh_k^\zZ(\aA)$ of a subcategory $\aA \subset \zZ$. Using the language of kernel and cokernel pairs (of subcategories in $\zZ$) we show that $\ngbh_k^\zZ(\aA)$ is abelian if $\aA \subset \zZ$ is a full abelian subcategory.

\vspace{0.3cm}
\subsection{Gabriel product and infinitesimal extensions}~\\

Let $\zZ$  be an exact category \cite{Quillen}, \cite{Buehl}, for example, an abelian one. We will use the terminoloy of {\em conflations, deflations and inflations} due B. Keller \cite{Kel4}. For strictly full subcategories $\aA$ and $\bB$ in $\zZ$, consider their \emph{Gabriel product}
$$
\aA \ast \bB,
$$
i.e. the full subcategory of $\zZ$ whose objects are $C \in \zZ$ which fit into  conflations $A \to C\to B$ with $A \in \aA$, $B\in \bB$.

\begin{LEM}\label{lem_Gab_prod_ass}
	The Gabriel product is associative, i.e. $(\aA \ast \bB) \ast \cC = \aA \ast(\bB \ast \cC)$.
\end{LEM}
\begin{proof}
	Given $Z\in (\aA \ast \bB) \ast \cC$, one can construct a commutative diagram with $A\in \aA$, $B\in \bB$, $C\in \cC$ and conflations in rows and columns \cite[Corollary 3.6]{Buehl}:
	\[
	\xymatrix{0 \ar[r] & C \ar[r]^{\Id} & C\\
	A \ar[r] \ar[u] & Z \ar[r] \ar[u] & Z_2 \ar[u] \\ A \ar[r] \ar[u]^{\Id} & Z_1 \ar[r] \ar[u] & B \ar[u]}
	\]
	which proves that $Z$ belongs to 
	$\aA \ast (\bB\ast \cC)$, and vice versa.
\end{proof}

Let $\aA \subset \zZ$ be a strictly full subcategory.
The \emph{$k$-infinitesimal extension} of $\aA$ in $\zZ$ is the full subcategory 
$$
\ngbh_{k}^{\zZ}(\aA):=\underbrace{\aA\ast \ldots \ast \aA}_{k+1}.
$$
We put $\ngbh_0^{\zZ}(\aA) = \aA$. 
Objects of $\ngbh_k^{\zZ}(\aA)$ are said to be  \emph{$k$-infinitesimal with respect to} $\aA$. We say that an object of $\zZ$ is \emph{infinitesimal with respect to $\aA$} if it is $k$-infinitesimal for some $k\in \mathbb{N}$.

We say that a strictly full subcategory $\aA\subset \zZ$ 
is \emph{extension closed} if, for any $A, A'\in \aA$, any extension of $A$ by $A'$ is an object of $\aA$, i.e. if $\ngbh^{\zZ}_1(\aA) \simeq \aA$.
We define the \emph{extension closure}  of $\aA$ in $\zZ$ as the full subcategory of objects infinitesimal with respect to $\aA$, and denote the closure by $\ngbh^{\zZ}(\aA): = \bigcup_k \ngbh^{\zZ}_k(\aA)$.  Clearly,  $\ngbh^{\zZ}(\aA)$ is extension closed.
When $\zZ$ is clear from the context, we drop $\zZ$ from the notation and write $\ngbh_k(\aA)$ or $\ngbh(\aA)$ .

We say that a strictly full subcategory $\aA$ in  $\zZ$ is \emph{extension generating}, or simply \emph{ext-generating}, if $\ngbh^{\zZ}(\aA) \simeq \zZ$. Analogously, a fully faithful functor $F\colon \aA\to \zZ$ is \emph{extension generating} (\emph{ext-generating}) if so is its essential image $F(\aA)\subset \zZ$.

\vspace{0.3cm}
\subsection{Kernel and cokernel pairs}~\\

Let now $\zZ$ be an abelian category and $\aA, \bB \subset \zZ$ strictly full subcategories. We say that \emph{$\aA$ has kernels relative to $\bB$} if, for any $A\in \aA$, $B\in \bB$ and $f\in \Hom_{\zZ}(A,B)$, the kernel of $f$ belongs to $\aA$. We say that \emph{$\bB$ has cokernels relative to $\aA$} if, for any $A\in \aA$, $B \in \bB$ and $f\in \Hom_{\zZ}(A,B)$, the cokernel of $f$ belongs to $\bB$.

We say that a pair $(\aA, \bB)$ of strictly full subcategories $\zZ$ is a \emph{kernel pair} if $\aA$ has kernels relative to $\bB$ (resp. a \emph{cokernel pair} if $\bB$ has cokernels relative to $\aA$).

\begin{LEM}\label{lem_AB_C}
	For strictly full subcategories $\aA , \bB , \cC$ in $\zZ$, the pair $(\aA \ast \bB, \cC)$ is 
	\begin{enumerate}
		\item a kernel pair if $(\aA, \cC)$ is a kernel and cokernel pair and $(\bB, \cC)$ is a kernel pair,
		\item a cokernel pair if $(\aA, \cC)$ and $(\bB, \cC)$ are cokernel pairs.
	\end{enumerate}
\end{LEM}
\begin{proof}
	Consider $C\in \cC$ and $Z\in \aA\ast \bB$ which fits into a conflation $A \xrightarrow{i} Z\xrightarrow{p} B$ with $A\in \aA$ and $B\in \bB$. A morphism $f\colon Z \to C$ yields a morphism of conflations:
	\[
	\xymatrix{C \ar[r]^{\Id} & C \ar[r] & 0\\
	A \ar[r]^i \ar[u]^{fi} & Z \ar[r]^p \ar[u]^f & B. \ar[u]^0}
	\]
	By the snake lemma \cite[Corollary 8.13]{Buehl}, we get an exact sequence:
	$$
	0 \to \ker fi \to \ker f \to B \to \coker fi \to \coker f \to 0.
	$$
	Conditions \emph{(1)} imply that $\ker fi \in \aA$, $\coker fi \in \cC$, and $\ker (B \to \coker fi)\in \bB$. Hence $\ker f\in \aA\ast \bB$.
	
	Conditions \emph{(2)} imply that $\coker fi\in \cC$, and $\coker f = \coker (B \to \coker fi) \in \cC$.
\end{proof}

\begin{LEM}\label{lem_C_AB}
		For strictly full subcategories $\aA , \bB , \cC$ in $\zZ$, the pair $(\cC, \aA \ast \bB)$ is 
	\begin{enumerate}
		\item  a kernel pair if $(\cC, \aA)$ and $(\cC, \bB)$ are kernel pairs,
		\item  a cokernel pair if $(\cC, \aA)$ is a cokernel pair and $(\cC, \bB)$ is a kernel and a cokernel pair.
	\end{enumerate}
\end{LEM}
\begin{proof}
	Consider $C\in \cC$ and $Z\in \aA\ast \bB$ which fits into a short exact sequence $0 \to A \xrightarrow{i} Z \xrightarrow{p} B \to 0$ with $A\in \aA$ and $B\in \bB$. A morphism $f\colon C\to Z$ yields a morphism of short exact sequences
	\[
	\xymatrix{A \ar[r]^i & Z \ar[r]^p & B\\
		0 \ar[r] \ar[u] & C \ar[r]^{\Id} \ar[u]^f & C. \ar[u]^{pf}}
	\]
	By the snake lemma we get an exact sequence:
	$$
	0 \to \ker f \to \ker pf \to A \to \coker f \to \coker pf \to 0.
	$$
	Conditions \emph{(1)} imply that $\ker pf \in \cC$, hence $\ker f = \ker(\ker pf \to A) \in \cC$.
	
	Conditions \emph{(2)} imply that $\coker pf \in \bB$, $\ker pf \in \cC$ and $\coker(\ker pf \to A)\in \aA$. It follows that $\coker f \in \aA \ast \bB$.
\end{proof}

\begin{PROP}\label{prop_ker_pair}
	The pair $(\aA \ast \bB, \aA \ast \bB)$ is 
	\begin{enumerate}
		\item a kernel pair if $(\aA, \bB)$ and $(\aA, \aA)$ are kernel and cokernel pairs and $(\bB, \bB)$ and $(\bB, \aA)$ are kernel pairs,
		\item a cokernel pair if $(\bB, \bB)$ and $(\aA, \bB)$ are kernel and cokernel pairs and $(\bB, \aA)$ and $(\aA, \aA)$ are cokernel pairs.
	\end{enumerate}
\end{PROP}
\begin{proof}
	By Lemma \ref{lem_C_AB}, $(\aA \ast \bB, \aA \ast \bB)$ is a kernel pair if $(\aA \ast \bB, \aA)$ and $(\aA \ast \bB, \bB)$ are. Statement \emph{(1)} follows from Lemma \ref{lem_AB_C}.
	
	By Lemma \ref{lem_AB_C}, $(\aA \ast \bB, \aA\ast \bB)$ is a cokernel pair if $(\aA \ast \bB, \aA)$ and $(\aA \ast \bB, \bB)$ are. Statement \emph{(2)} follows from Lemma \ref{lem_C_AB}.
\end{proof}

Let $\zZ$ be an abelian category. We say that $\aA\subset \zZ$ is a \emph{full abelian subcategory} if
$\aA$ is abelian and the embedding functor $\aA\to \zZ$ is exact and fully faithful, i.e. if $(\aA, \aA)$ is both a kernel and a cokernel pair. Note that an abelian full subcategory $\aA \subset \zZ$, i.e. a full subcategory that is abelian, is not necessarily a full abelian subcategory as kernels and cokernels in $\aA$ and $\zZ$ may differ. 

\begin{COR}
	$\aA \ast \bB$ is a full abelian subcategry if $\aA$ and $\bB$ are full abelian subcategories and $(\aA, \bB)$ and $(\bB, \aA)$ are kernel and cokernel pairs.
\end{COR}
 
\vspace{0.3cm}

\vspace{0.3cm}
\subsection{Infinitesimal extensions  of full abelian subcategories are abelian}~\\

\begin{PROP}\label{prop_k-nghb_is_abelian}
	The $k$-infinitesimal extension $\ngbh_k^\zZ(\aA)$ of a full abelian subcategory $\aA \subset \zZ$ is abelian.
\end{PROP}
\begin{proof}
	We check that $\ngbh_k^\zZ(\aA)$ is closed under kernels and cokernels. 
	
	For a morphism $f\colon Z_1 \to Z_2$ in $\ngbh_k^{\zZ}(\aA)$, let $l, n\in \mathbb{N}$ be such that $Z_1\in \ngbh_{l}^{\zZ}(\aA)$, $Z_2\in \ngbh_{n}^{\zZ}(\aA)$. We check that $\ker f \in \ngbh_{l}^{\zZ}(\aA)$, $\coker f \in \ngbh_{n}^{\zZ}(\aA)$, i.e. that $(\ngbh_{l}(\aA), \ngbh_{n}(\aA))$ is a kernel and a cokernel pair.
	
	First, we show by induction on $n$ that $(\aA, \ngbh_{n}(\aA))$ is a kernel and a cokernel pair. Since $\aA \subset \zZ$ is a full abelian subcategory, $(\aA, \aA)$ is a kernel and a cokernel pair. By Lemma \ref{lem_Gab_prod_ass}, $\ngbh_n(\aA) \simeq \ngbh_{n-1}(\aA) \ast \aA$. Then,  Lemma \ref{lem_C_AB} and the inductive hypothesis imply that $(\aA, \ngbh_n(\aA))$ is a kernel and a cokernel pair.

	Next, for fixed  $n$ we show by induction on $l$ that $(\ngbh_l(\aA), \ngbh_n(\aA))$ is a kernel and a cokernel pair. The case $l=0$, i.e. $(\aA, \ngbh_n(\aA))$ is discussed above. By Lemma \ref{lem_Gab_prod_ass}, 
	$\ngbh_{l}(\aA) \simeq \ngbh_{l-1}(\aA) \ast \aA$. The statement follows from Lemma \ref{lem_AB_C}.
\end{proof}
 
Since an ascending union of abelian subcategories is abelian, Proposition \ref{prop_k-nghb_is_abelian} implies
 \begin{COR}\label{cor_nghb_is_abelian}
 	The extension closure $\ngbh^{\zZ}(\aA)$ of a full abelian subcategory $\aA\subset \zZ$ is abelian.
 \end{COR}

\section{Non-commutative deformations of an exact functor}\label{sec_nc_def_exact_fc}~\\

We introduce the category $\wt{\textrm{Fin}}_{\aA_0}$ of abelian categories which are extension generated by an abelian subcategory equivalent to $\aA_0$. Category $\wt{\textrm{Fin}}_{\aA_0}$ has a full subcategory $\Fin_{\aA_0}$ whose objects are categories equivalent to $\ngbh_k(\aA_0)$, for some $k\in \mathbb{N}$. 
We define the functor $\ncDef_\zeta\colon \Fin_{\aA_0}^{\opp} \to \Sets$ of non-commutative deformations of a given exact functor $\zeta\colon \aA_0 \to \zZ$. We prove that $\wt{\Fin}_{\aA_0}$ is a full subcategory of the category of ind-objects over $\Fin_{\aA_0}$ and that 
 the category $\ngbh^\zZ(\aA_0)\in \wt{\Fin}_{\aA_0}$ ind-represents $\ncDef_\zeta$ if the functor $\zeta$ is fully faithful. We give an example of an exact functor $\xi\colon\modu \kk^3 \to \zZ$ for which $\ncDef_\xi$ is not ind-representable.

\vspace{0.3cm}
\subsection{Definition of the deformation functor}\label{ssec_def_of_def}~\\

For an abelian category $\aA_0$, we denote by $\wt{\Fin}_{\aA_0}$ the category of abelian categories extension generated by $\aA_0$. Objects in $\wt{\Fin}_{\aA_0}$ are pairs $(\aA, Q_{\aA})$ of an abelian category $\aA$ and a fully faithful extension generating exact functor $Q_{\aA} \colon \aA_0\to \aA$. For $(\aA, Q_{\aA}), (\bB, Q_{\bB})$ in $\wt{\Fin}_{\aA_0}$,
\begin{align*}
\Hom_{\wt{\Fin}_{\aA_0}}((\aA, Q_{\aA}), (\bB, Q_{\bB})) = \{(F, \f \colon F\circ Q_{\aA} \xrightarrow{\simeq} Q_{\bB})| F\colon \aA\to \bB \emph{ is exact}\}/\sim
\end{align*}
where $(F,\f) \sim (F', \f')$ if there exists $\psi \colon F\xrightarrow{\simeq} F'$ such that $\f'\circ \psi_{Q_\aA}= \f$. We often omit $Q_{\aA}$ and write $\aA\in \wt{\Fin}_{\aA_0}$.

We denote by $\Fin_{\aA_0}$ the full subcategory of $\wt{\Fin}_{\aA_0}$ of \emph{$\aA_0$-finite categories}, i.e. those $\aA\in \wt{\Fin}_{\aA_0}$ for which there exist $k\in \mathbb{N}$ such that $ \ngbh^k(\aA_0) \xrightarrow{\simeq} \aA$.

Note that $\wt{\Df_n} \subset \wt{\Fin}_{\modu \kk^n}$, $\Df_n \subset \Fin_{\modu \kk^n}$ are full subcategories whose objects are $\Ext^1$-finite categories (see
Theorem \ref{thm_Df_as_bound_on_Loewy}).

Given an exact functor $\zeta \colon \aA_0 \to \zZ$ we define a functor $\ncDef_\zeta \colon \Fin_{\aA_0}^{\opp} \to \Sets$ via:
\begin{equation}\label{eqtn_ncDef_cat_gen}
\ncDef_{\zeta}(\aA, Q_{\aA}) = \{(T,\tau\colon T\circ Q_{\aA}\xrightarrow{\simeq} \zeta )\,|\, T\colon \aA\to \zZ \textrm{ is exact}\}/\sim
\end{equation}
where $(T,\tau) \sim (T',\tau')$ if there exists $\psi \colon T\xrightarrow{\simeq}T'$ such that $\tau'\circ \psi_{ Q_\aA} = \tau$. 

For $(F,\f)\in \Hom_{\Df_n}(\aA, \bB)$, we 
define $\ncDef_\zeta(F,\f) \colon \ncDef_{\zeta}(\bB) \to \ncDef_{\zeta}(\aA)$ via:
\begin{equation}
\ncDef_\zeta(F,\f)(T,\tau) = (T\circ F, T\circ F \circ Q_\aA \xrightarrow{T\f} T\circ Q_\bB \xrightarrow{\tau} \zeta).
\end{equation}

\vspace{0.3cm}
\subsection{$\wt{\Fin}_{\aA_0}$ as a subcategory of ind-objects over $\Fin_{\aA_0}$}\label{ssec_Fin_as_ind-obj}~\\

Let $\cC$ be a category and $\cC^{\wedge}$ the category of functors $\cC^{\textrm{op}} \to \textrm{Sets}$.
Consider the Yoneda embedding $h\colon \cC \to \cC^{\wedge}$, $C\mapsto h_C(-) \colon = \Hom_{\cC}(-,C)$. An \emph{ind-object} over $\cC$ is an object of $\cC^\wedge$ isomorphic to $\varinjlim_I h \circ \alpha$, for some functor $\alpha \colon I \to \cC$ with $I$ filtrant and small. 
We denote by $\textrm{Ind}(\cC)$ the full subcategory of $\cC^\wedge$ with ind-objects as objects.  

We say that a functor $\Upsilon \in \cC^\wedge$ is \emph{ind-representable} if it is isomorphic to an ind-object.

\begin{LEM}\label{lem_homs_to_Fin_n}
	For $\aA\in \Fin_{\aA_0}$, $\bB\in \wt{\Fin}_{\aA_0}$ and $T\in \Hom_{\wt{\Fin}_{\aA_0}}(\aA, \bB)$ there exists $l\in \mathbb{N}$ such that $T(\aA) \subset \ngbh_l^{\bB}(\aA_0)$ and functor $T$ factors via the inclusion $i_l\colon \ngbh_l^{\bB}(\aA_0) \to \bB$. Therefore,
	$$
	\Hom_{\wt{\Fin}_{\aA_0}}(\aA, \bB) = \varinjlim \Hom_{\Fin_{\aA_0}}(\aA, \ngbh_l^{\bB}(\aA_0)).
	$$
\end{LEM}
\begin{proof}
	Any functor $F \in \Hom_{\wt{\Fin}_{\aA_0}}(\aA, \bB)$ restricts to a a functor $F_l \colon \ngbh_l^{\aA}(\aA_0) \to \ngbh_l^{\bB}(\aA_0)$. Indeed, the isomorphism $F\circ Q_{\aA} \simeq Q_{\bB}$ implies that functor $F$ induces an equivalence of the subcategory $\aA_0 \subset \aA$ with the subcategory $\aA_0 \subset \bB$. Exactness of $F$ implies that $F|_{\ngbh_l^{\aA}(\aA_0)}$ is a functor $\ngbh_l^{\aA}(\aA_0) \to \ngbh_l^{\bB}(\aA_0)$.
		
	By definition an $\aA_0$-finite category $\aA \in \Fin_{\aA_0}$ is equivalent to $\ngbh_l^{\aA}(\aA_0)$, for some $l\in \mathbb{N}$. It follows that, for $\aA \in \Fin_{\aA_0}$ and $\bB\in \wt{\Fin}_{\aA_0}$,  a functor $F \in \Hom_{\wt{\Fin}_{\aA_0}}(\aA, \bB)$ decomposes as $F\colon \aA \simeq \ngbh_l^{\aA}(\aA_0)\xrightarrow{F_l} \ngbh_l^{\bB}(\aA_0) \xrightarrow{i_l}\bB$.
\end{proof}

\begin{THM}\label{thm_ind_object_general}
	 The functor $h\colon \wt{\Fin}_{\aA_0} \to {\Fin_{\aA_0}}^\wedge$, $h_{\aA}(-) = \Hom_{\wt{\Fin}_{\aA_0}}(-,\aA)$, for $\aA \in  \wt{\Fin}_{\aA_0}$, induces an equivalence of $\wt{\Fin}_{\aA_0}$ with a full subcategory of $\Ind(\Fin_{\aA_0})$.
\end{THM}
\begin{proof}
	Lemma \ref{lem_homs_to_Fin_n} implies that, for any $\aA \in \wt{\Fin}_{\aA_0}$, functor $h_{\aA}(-)$ is an ind-object over $\Fin_{\aA_0}$.
	It remains to check that functor $h $ is fully faithful.
	
	By \cite[formula (2.6.4)]{KasSch2}, $\Hom_{\Fin_{\aA_0}^\wedge}(h_{\aA}, h_{\bB}) \simeq \varprojlim_i \varinjlim_j \Hom_{\Fin_{\aA_0}}(\ngbh_i^{\aA}(\aA_0), \ngbh_j^{\bB}(\aA_0))$. 
	Lemma \ref{lem_homs_to_Fin_n} implies that $\Hom_{{\Fin_{\aA_0}}^\wedge}(h_{\aA}, h_{\bB}) \simeq \varprojlim \Hom_{\wt{\Fin}_{\aA_0}}(\ngbh_i^{\aA}(\aA_0), \bB)$. Under this isomorphism, the morphism of Hom-spaces $h_{\aA, \bB} \colon \Hom_{\wt{\Fin}_{\aA_0}}(\aA, \bB) \to \varprojlim \Hom_{\wt{\Fin}_{\aA_0}}(\ngbh_i^{\aA}(\aA_0),\bB)$ induced by $h$ maps $F \in \Hom_{\wt{\Fin}_{\aA_0}}(\aA,\bB)$ to the family $\{F|_{\ngbh_i(\aA_0)}\in \Hom_{\wt{\Fin}_{\aA_0}}(\ngbh_i^{\aA}(\aA_0), \bB)\}_{i}$.

	As any object  and any morphism in $\aA$ belongs to the subcategory $\ngbh_i^{\aA}(\aA_0)$, for some $i$, morphism $h_{\aA, \bB}$ is an isomorphism, i.e. functor $h$ is fully faithful.
\end{proof}

Analogous arguments show
\begin{THM}\label{thm_ind_object}
	The functor $h\colon \wt{\Df}_{n} \to {\Df_{n}}^\wedge$, $h_{\aA}(-) = \Hom_{\wt{\Df}_{n}}(-,\aA)$ induces an equivalence of $\wt{\Df}_{n}$ with a full subcategory of $\Ind(\Df_{n})$.
\end{THM}

\vspace{0.3cm}
\subsection{Ind-representability of nc deformations of fully faithful functors}\label{ssec_ind_rep_ff}~\\

Let $\zeta\colon \aA_0 \to \zZ$ be a fully faithful exact functor of abelian categories. By Corollary \ref{cor_nghb_is_abelian}, $\ngbh^\zZ(\aA_0)$ is an abelian category. Denote by $\iota \colon \ngbh^\zZ(\aA_0) \to \zZ$ the inclusion functor. Then $\zeta$ decomposes as $\iota \circ \zeta_0$, for $\zeta_0 \colon \aA_0 \to \ngbh^{\zZ}(\aA_0)$, and $(\ngbh^{\zZ}(\aA_0), \zeta_0)$ is an object of $\wt{\Fin}_{\aA_0}$.

\begin{THM}\label{thm_ind_rep_fully_faith}
	Consider a fully faithful exact functor $\zeta \colon \aA_0 \to \zZ$ of abelian categories. Then $(\ngbh^{\zZ}(\aA_0), \zeta_0)$ ind-represents the functor $\ncDef_\zeta\colon \Fin_{\aA_0}^{\opp} \to \Sets$.
\end{THM}
\begin{proof}
	The inclusion $\iota$ induces a morphism of functors $\Phi\colon h_{\ngbh(\aA_0)} \to \ncDef_\zeta$ which maps $(T, \tau \colon T \circ Q_{\aA} \xrightarrow{\simeq} \zeta_0)$ in $\Hom_{\wt{\Fin}_{\aA_0}}(\aA, \ngbh(\aA_0))$ to $(\iota \circ T, \iota \circ T \circ Q_{\aA} \xrightarrow{\iota \tau} \iota \circ \zeta_0\simeq \zeta) \in \ncDef_{\zeta}(\aA)$.
	
	As any exact functor $F\colon \aA\to \zZ$ such that $F\circ Q_{\aA} \simeq \zeta$ is uniquely decomposed as $\iota \circ \ol{F}$, for $\ol{F} \colon \aA\to \ngbh(\aA_0)$, the natural transformation $\Phi$ is an isomorphism.
\end{proof}

In the case when $\aA_0 \simeq \modu \kk^n$ and category $\zZ$ is $\kk$-linear, an (exact) functor $\zeta \colon \modu \kk^n \to \zZ$ is determined by a $\kk^n$-object in $\zZ$, see  Proposition \ref{prop_functor_given_by_tensor}, i.e. by an $n$-tuple $\sigma =(Z_1,\ldots,Z_n)$ of objects in $\zZ$. Functor $\zeta$ is fully faithful if and only if $\sigma$ is a \emph{simple collection} (\emph{cf.} \cite{Kaw5}), i.e. if and only if
$$
\dim_{\kk}\Hom_\zZ(Z_i,Z_j) =\delta_{ij}.
$$
Theorem \ref{thm_ind_rep_fully_faith} implies
\begin{COR}\label{cor_simpl_coll_gen}
	Let $\sigma =(Z_1,\ldots,Z_n)$ be a simple collection in a $\kk$-linear abelian category $\zZ$ and $\zeta \colon \modu \kk^n \to \zZ$ an exact functor determined by it. If $\sigma$ is a simple collection then $\ncDef_{\zeta} \colon \Fin_{\modu \kk^n}^{\opp} \to \Sets$ is ind-represented by the extension closure $\ngbh^\zZ(\sigma) \in \wt{\Fin}_{\modu \kk^n}$.
\end{COR}

The ind-representability of the deformation functor $\ncDef_\zeta \colon \Df_n^{\opp} \to \Sets$ (\ref{eqtn_ncDef_cat}) requires the deformed collection $\sigma$ to satisfy a finiteness condition:

\begin{THM}\label{thm_pro-represent_for_simple_coll}
	Let $\sigma = (Z_1,\ldots, Z_n)$ be a simple collection in a $\kk$-linear abelian category $\zZ$ and $\zeta \colon \modu \kk^n \to \zZ$ an exact functor determined by it. If $\sigma$ is a simple collection and $\dim \Ext^1_{\zZ}(\bigoplus_i Z_i, \bigoplus_i Z_i)$ is finite, then the category $\ngbh^{\zZ}(\sigma) \in \wt{\Df}_n$ ind-represents the functor $\ncDef_{\zeta} \colon {\Df_n}^{\opp} \to \Sets$.
\end{THM}
\begin{proof}
	Any object in $\ngbh^\zZ(\sigma)$ admits a finite filtration with graded factors in $\sigma$. Hence, if $\Ext^1_{\zZ}(\bigoplus Z_i, \bigoplus Z_i)$ is finite dimensional then $\ngbh^\zZ(\sigma)$ is $\Ext^1$-finite, i.e.  $\ngbh^\zZ(\sigma)\in \wt{\Df}_n$ (see  Theorem \ref{thm_Df_as_bound_on_Loewy}). The statement follows from Corollary \ref{cor_simpl_coll_gen}.
\end{proof}

\vspace{0.3cm}
\subsection{A counterexample to ind-representability of the deformation functor}\label{ssec_obst_to_rep}~\\

We give an example of a not fully faithful functor $\zeta \colon \modu \kk^{\oplus 3} \to \zZ$ such that $\ncDef_{\zeta}$ is not ind-representable.

Consider category $\bB$ of right modules over the path algebra $B$ of the quiver
\begin{equation}\label{eqtn_quiver_B}
\xymatrix{1 & 2 \ar[l] & 3 \ar[l]}
\end{equation}
and $Q_B={q_B}_*$ is given by a fixed homomorphism $q_B \colon B\to \kk^{\oplus 3}$ inducing an isomorphism $B/\textrm{rad}B \xrightarrow{\simeq} \kk^{\oplus 3}$. 

Let $\aA$ be the category of right modules over the path algebra $A$ of the quiver
\begin{equation}\label{eqtn_quiver_A}
	\xymatrix{1 & 2  & 3. \ar[l]}
\end{equation}
Let $\alpha \colon B\to A$ be the algebra homomorphism that kills the arrow $1 \leftarrow 2$. It induces the inclusion of the subcategory $\alpha_*\colon \aA\to \bB$. Let $q_A \colon A\to \kk^{\oplus 3}$ be such that $q_A \circ \alpha = q_B$ and let $\gamma$ be an isomorphism of functors $\alpha_*{q_A}_* \xrightarrow{\simeq}{q_B}_*$ induced by this equality. Then $(\alpha_*, \gamma)$ is an element in $\Hom_{\Df_3}((\aA, Q_{\aA}), (\bB, Q_{\bB}))$, where $Q_{\aA}= {q_A}_*$.

An $A$-object in a $\kk$-linear abelian category $\wW$ is a functor from the path category of the quiver (\ref{eqtn_quiver_A}), i.e. the data of $\rho_A = (W'_1, W'_2 \xleftarrow{\sigma} W'_3)$. Similarly, a $B$-object is $\rho_B = (W_1 \xleftarrow{\tau_1} W_2 \xleftarrow{\tau_2} W_3)$. The corresponding right exact functor $T_{\rho_A}$ (see Proposition \ref{prop_left_adj_to_Hom}), respectively $T_{\rho_B}$, maps the indecomposable projective covers $P_i$ of simple module $S_i$ corresponding to the $i$'th vertex of (\ref{eqtn_quiver_A}), respectively (\ref{eqtn_quiver_B}), to $W'_i$, respectively to $W_i$, and the non-trivial morphisms $P_i \to P_j$ to given morphisms $W'_i \to W'_j$, respectively $W_i \to W_j$.

Indecomposable projective covers of simple objects $S_2^A$ and $S_3^A$ in $\aA$ become, under $\alpha_*$, projective in $\bB$ too. The simple object $S_1^A$ is projective in $\aA$, while in $\bB$ a projective resolution of $S_1^B$ is: 
\begin{equation}\label{eqtn_res_of_S1}
0 \to P_2^B \to P_1^B \to S_1^B \to 0.
\end{equation} 
\begin{LEM}
	An $A$-object $\rho_A$ is flat if and only if $\sigma$ is a monomorphism. A $B$-object $\rho_B$ is flat if and only if $\tau_1$ and $\tau_2$ are monomorphisms.
\end{LEM}
\begin{proof}
	We prove the first statement. The proof of the second one is analogous.
	
	Category $\aA$ has 4 isomorphism classes of indecomposable objects: simple $A$-modules $S_1^A$, $S_2^A$, $S_3^A$ and the projective cover $P_2^A$ of $S_2$ which fits into a short exact sequence 
	\begin{equation}\label{eqtn_res_S_2}
	0 \to S_3^A \xrightarrow{i} P_2^A \xrightarrow{\pi} S_2^A\to 0
	\end{equation}
	Functor $T_{\rho_A}\colon \aA\to \wW$ is defined as $T_{\rho_A}(S_1^A)= W'_1$, $T_{\rho_A}(S_2^A) = \textrm{coker } \sigma$, $T_{\rho_A}(S_3^A)= W'_3$, $T_{\rho_A}(P_2^A) = W'_2$, $T_{\rho_A}(i) = \sigma$ and $T_{\rho_A}(\pi)$ is the canonical map $W'_2 \to \textrm{coker }\sigma$.
	
	If $T_{\rho_A}$ is exact then it maps (\ref{eqtn_res_S_2}) to a short exact sequence. In particular, $\sigma$ is a monomorphism. If, on the other hand, $\sigma$ is a monomorphism, then $L^1T_{\rho_A}$ vanishes on all simple objects in $\aA$. Hence, $L^1T_{\rho_A} =0$, which implies that $T_{\rho_A}$ is exact.
\end{proof}

Consider a $\kk$-linear abelian category $\wW$ and an exact functor $\eta \colon \kk^{\oplus 3}{\rm -mod }\to \wW$. Let $(S, \sigma)$ be an element of $\ncDef_{\eta}(\aA, Q_{\aA})$.
Assume that $S$ is given by a flat  $A$-object $(W_1, W_2 \xleftarrow{s} W_3)$. Consider
\begin{equation}\label{eqtn_theta_W}
	\theta_{\wW} \colon \Ext^1_{\wW}(W_1, W_2) \to \ncDef_{\eta}(\bB, Q_{\bB})\end{equation}
defined as follows: 
for an  element $\xi\in \Ext^1_{\wW}(W_1, W_2)$, corresponding to an exact sequence 
$$
0 \to W_2\xrightarrow{t} \wt{W_1} \to W_1\to 0,
$$
 $\theta_{\wW}(\xi) =(G_\xi, \rho_\xi)$, where
$G_{\xi}\colon \bB \to \wW$ is the functor given by the (flat) $B$-object $(\wt{W_1} \xleftarrow{t} W_2 \xleftarrow{s} W_3)$ and $\rho_\xi\colon G_\zeta \circ Q_{\bB} \to \eta$ is induced by $\sigma\colon S \circ Q_{\aA} \to \eta$ and the canonical isomorphism $G_{\xi}|_{\aA}\xrightarrow{\simeq} S$.

\begin{REM}\label{rem_mu_id_on_simple}
	Recall that $(G_{\xi_1}, \rho_{\xi_1})$ and $(G_{\xi_2}, \rho_{\xi_2})$ define the same element in $\ncDef_{\eta}(\bB)$ if and only if there exists an isomorphism $\psi \colon G_{\xi_1} \xrightarrow{\simeq} G_{\xi_2}$ such that $\rho_{\xi_2}\circ (\psi{Q_{\bB}}) = \rho_{\xi_1}$. Note that, as isomorphisms $\rho_{\xi_1}$ and $\rho_{\xi_2}$ are induced by $\sigma$, the last equality holds if and only if $\psi$ induces the identity natural transformation on the essential image of $\eta$. 
\end{REM}

Note that if $(S, \sigma) \in \Hom_{\wt{Df}_3}((\aA, Q_{\aA}), (\wW, \eta))$ then $\theta_{\wW}$ takes values in the set $\Hom_{\wt{\Df}_3}((\bB, Q_{\bB}),(\wW, \eta))$.
\begin{LEM}\label{lem_theta_inj} Assume that $(\wW, \eta)\in\wt{\Df}_3$. If $(S, \sigma) \in \Hom_{\wt{Df}_3}((\aA, Q_{\aA}), (\wW, \eta))$, then $\theta_{\wW}$ is injective.
\end{LEM}
\begin{proof}
	By assumption $(\wW, \eta)$ is an object in $\wt{\Df}_3$. In particular, objects $W_1$, the cokernel of $s$, and $W_3$ are pairwise non-isomorphic simple objects in $\wW$.

	Consider $\xi_1, \xi_2 \in \Ext^1_{\wW}(W_1, W_2)$ and corresponding exact sequences
	\begin{align*}
		& 0 \to W_2 \xrightarrow{t_1} \wt{W_1^1} \to W_1 \to 0,& & 0 \to W_2 \xrightarrow{t_2} \wt{W_1^2} \to W_1 \to 0.&
	\end{align*} 
Assume that $\psi \colon G_{\xi_1} \xrightarrow{\simeq} G_{\xi_2}$ is an isomorphism such that $\rho_{\xi_2} \circ \psi_{Q_{\bB}} = \rho_{\xi_1}$. By Remark \ref{rem_mu_id_on_simple}, $\psi$ is  the identity natural transformation on simple objects in $\wW$. Hence, $\psi_2 \colon W_2\to W_2$ is a morphism that induces the identity morphism on the quotient $W_2/W_3$.

We want to show that $\psi_2 = \Id_{W_2}$.
Assume first that $W_2$ is a non-trivial extension of simple objects in $\wW$. Then the endomorphism algebra of $W_2$ is isomorphic to $\kk$.  As $\la \Id_{W_2}$ induces $\la \Id$ on the maximal semisimple quotient $W_2/W_3$ of $W_2$, we conclude that $\psi= \Id_{W_2}$ is the identity morphism. If, on the other hand,  $W_2\simeq W_3\oplus W_2/W_3$ then $\psi_2$ is the direct sum of the isomorphisms of simple direct summands, hence it is the identity morphism too.

It follows that isomorphism $\psi$ yields a morphism of short exact sequences
\[
\xymatrix{0 \ar[r] & W_2 \ar[r]^{t_2} & \wt{W_1^2} \ar[r] & W_1 \ar[r] & 0 \\
	0 \ar[r] & W_2 \ar[r]^{t_1} \ar[u]^{\Id_{W_2}}& \wt{W_1^2} \ar[r] \ar[u]^{\psi_1} & W_1 \ar[r] \ar[u]^{\Id_{W_1}} & 0 }
\]
Hence, classes of $\xi_1$ and $\xi_2$ are equal in $\Ext^1_{\wW}(W_1,W_2)$.
\end{proof}

 Consider the category $\zZ$ of right finitely dimensional modules over the quiver
\begin{equation}\label{eqtn_quiver_Z}
	\xymatrix{1 & 3 \ar[l] \ar@<1ex>[r]& 2. \ar@<1ex>[l]}
\end{equation}
We denote by $L_1$, $L_2$, $L_3$ simple objects in $\zZ$ and consider a non-trivial extension
$$
0 \to L_2 \to M \xrightarrow{a} L_3 \to 0.
$$
Consider functor $\zeta\colon \modu \kk^{\oplus 3} \to \zZ $  given by the $\kk^{\oplus 3}$-object:
\begin{equation}\label{eqtn_Zobject}
Z_\bcdot=(L_1, M, L_3).
\end{equation} 
Functor $\zeta$ is not fully faithful as $a$ is a generator for $\Hom_{\zZ}(M,L_3)\simeq \kk$. The following extensions groups between $L_1$, $M$ and $L_3$ are non-zero:
\begin{align*}
	&\Ext^1_{\zZ}(L_1, M)\simeq \kk,& &\Ext^1_{\zZ}(M, L_3)\simeq \kk,& & \Ext^1_{\zZ}(L_1, L_3)\simeq \kk.&
\end{align*} 
Moreover, the composition 
\begin{align*}
	\Hom_{\zZ}(M, L_3) \otimes_{\kk} \Ext^1_{\zZ}(L_1, M) \to \Ext^1_{\zZ}(L_1, L_3)
\end{align*}
is an isomorphism.
Category $\zZ$ is hereditary, hence all higher Ext-groups vanish.

Consider
\begin{align*}
	0 \to L_3 \xrightarrow{b} N \xrightarrow{c} M \to 0 
\end{align*}
corresponding to a non-zero element in $\Ext^1_{\zZ}(M, L_3)$. Then
\begin{align*}
	&\End_{\zZ}(N) = \textrm{Span}_{\kk}(\Id_N, bac),& &(bac)^2 = 0,&
\end{align*}
i.e. $\End_{\zZ}(N)$ is the algebra of dual numbers. Moreover, the canonical action
\begin{align*}
	\End_{\zZ}(N) \otimes_{\kk} \Ext^1_{\zZ}(L_1, N)\to \Ext^1_{\zZ}(L_1, N)
\end{align*} 
gives $\Ext^1_{\zZ}(L_1, N)$ the structure of  a rank one free module over $\End_{\zZ}(N)$.

To prove that $\ncDef_{\zeta}$ is not ind-representable we argue by contradiction and assume that there exists $(\cC, Q_{\cC}) \in \wt{\Df}_3$ and 
$$
(T\colon \cC\to \zZ, \tau\colon T\circ Q_{\cC} \xrightarrow{\simeq} \zeta) 
$$
which induces a functorial isomorphism 
\begin{equation}\label{eqtn_iso_Psi}
\Psi \colon \Hom_{\wt{\Df}_3}(-,(\cC, Q_\cC)) \xrightarrow{\simeq}\ncDef_{\zeta}(-).
\end{equation}

 We consider functor $U\colon \aA\to \zZ$ given by the $A$-object $(L_1, N \xleftarrow{b} L_3)$. Let	$\upsilon \colon U\circ Q_{\aA} \to \zeta$  be the tautological isomorphism. 
 Then $(U, \upsilon) \in \ncDef_{\zeta}(\aA, Q_{\aA})$. As we assume that \eqref{eqtn_iso_Psi} is an isomorphism, there exists $(S, \sigma) \in \Hom_{\wt{\Df}_3}((\aA, Q_{\aA}), (\cC, Q_{\cC}))$ such that $\Psi(S, \sigma) = (T, \tau) \circ (S, \sigma) = (U, \upsilon)$. Let $C_1= S(S_1^A)$ and $C_2=S(P_2^A)$.
 
 Consider the commutative diagram
 \[
 \xymatrix{\Ext^1_{\zZ}(L_1, N) \ar[rr]^{\theta_{\zZ}} &&\ncDef_{\zeta}(\bB, Q_{\bB})\\
 	\Ext^1_{\cC}(C_1,C_2) \ar@{->>}[u]^{T\circ(-)} \ar@{^(->}[rr]^{\theta_{\cC}} & &\Hom_{\wt{\Df}_3}((\bB, Q_{\bB}), (\cC, Q_{\cC})) \ar[u]^{\Psi}_{\simeq} }
 \]
 where $\Psi$ is as in \eqref{eqtn_iso_Psi} and $\theta_{\zZ}$, $\theta_{\cC}$ as in \eqref{eqtn_theta_W}.
 
 The composition with $T$ yields a surjective map by Lemma \ref{lem_Ext_surj} below. Map $\theta_{\cC}$ is injective by Lemma \ref{lem_theta_inj}. As we assume that $\Psi$ is an isomorphism, then $\theta_{\zZ}$ must be injective too. In Lemma \ref{lem_theta_xi} we will show that $\theta_{\zZ}(\xi) = \theta_{\zZ}((\Id + bac)\xi)$, for any element $\xi\in \Ext^1_{\zZ}(L_1, N)$. As we can choose $\xi$ such that $\xi \neq (\Id + bac)\xi$, this contradicts the injectivity of $\theta_{\zZ}$. Thus, we have proven the following theorem modulo Lemmas \ref{lem_theta_xi} and \ref{lem_Ext_surj}.

\begin{THM}
Functor $\ncDef_\zeta$, for $\zeta$ as in \eqref{eqtn_Zobject}, is not ind-representable.
\end{THM}

\begin{LEM}\label{lem_theta_xi}
Consider $\theta_{\zZ}$ as above. Then, for any class $\xi$ in $\Ext^1_{\zZ}(L_1, N)$,  $\theta_{\zZ}(\xi) \sim \theta_{\zZ}((\mu \Id+\la bac)\xi)$ as elements of $\ncDef_{\zeta}(\bB, Q_{\bB})$ if and only if $\mu =1$.
\end{LEM} 
\begin{proof}
	First, assume that $\mu=1$. 
	By the definition of the action of $\End_{\zZ}(N)$ on $\Ext^1_{\zZ}(L_1,N)$ we have commutative diagram
	\[
	\xymatrix{0 \ar[r]&  N \ar[r]^{\wt{d}} &  \wt{N} \ar[r] & L_1 \ar[r] &  0 \\
	0 \ar[r] &  N \ar[r]^{\ol{d}} \ar[u]^{\Id + \la bac} &  \ol{N} \ar[r]\ar[u]^{\psi} & L_1 \ar[u]^{\Id}  \ar[r]&  0}
	\]
	where the lower row is the extension corresponding to $\xi$ and the upper row the extension corresponding to $(\Id + \la bac) \xi$. Moreover, $(\Id + \la bac)$ fits into a commutative diagram
	\[
	\xymatrix{0 \ar[r] & L_3 \ar[r]^b & N \ar[r]^c & M \ar[r] & 0 \\
	0 \ar[r] & L_3 \ar[r]^b \ar[u]^{\Id} & N \ar[r]^c \ar[u]^{\Id + \la bac} & M \ar[u]^{\Id} \ar[r]& 0}
	\]
	It follows that 
	\[
	\xymatrix{L_3 \ar[r]^b &  N \ar[r]^{\wt{d}} &  \wt{N} \\
	L_3 \ar[r]^b \ar[u]^{\Id} & N\ar[r]^{\ol{d}} \ar[u]^{\Id + \la bac} & \ol{N} \ar[u]^{\psi}}
	\]
	is an isomorphism of $B$-objects such that the induced isomorphism of functors $\theta_{\zZ}((\Id +\la bac)\xi) \to \theta_{\zZ}(\xi)$ is the identity on the image of $\zeta$. As the induced natural transformation on the essential image of $\zeta$ is the identity, $\theta_{\zZ}((\Id + \la bac)\xi)\sim  \theta_{\zZ}(\xi)$ as elements of $\ncDef_{\zeta}(\bB, Q_{\bB})$ (see Remark \ref{rem_mu_id_on_simple}). 
	
	Assume now that $\theta_{\zZ}(\xi) \sim \theta_{\zZ}((\mu \Id+\la bac)\xi)$, i.e. there exists an isomorphism $\psi$ of the corresponding functors $\bB \to \zZ$ which induces the identity natural transformation on the essential image of $\zeta$ (see Remark \ref{rem_mu_id_on_simple}). Then $\psi$ induces a morphism of short exact sequences
	\[
	\xymatrix{0 \ar[r] & N \ar[r] & \wt{N} \ar[r] & L_1 \ar[r] & 0\\
	0 \ar[r] & N \ar[r] \ar[u]^{\psi_N} & \ol{N} \ar[r] \ar[u]^{\psi_{\ol{N}}} & L_1 \ar[r] \ar[u]^{\Id} & 0,}
	\]
	where the first row is the extension corresponding to $(\mu \Id + \la bac) \xi$ and the second one to $\xi$. In particular, $\psi_{N} \xi = (\mu \Id + \la bac) \xi$. As $\Ext^1_{\zZ}(L_1, N)$ is a rank one free module over  $\End_{\zZ}(N)$, it implies that $\psi_N = \mu \Id + \la bac\in \End_{\zZ}(N)$. Then, as $M$ is the quotient of $N$ by $L_3$, the isomorphism $\psi_M\colon M \to M$ is $\mu \Id$. As we assumed that $\psi$ is the identity on the essential image of $\zeta$, we conclude that $\mu =1$.
\end{proof}

\begin{LEM}\label{lem_Ext_surj}
	Assume that $\Psi$ as in (\ref{eqtn_iso_Psi}) is an isomorphism. Then $T\circ (-) \colon \Ext^1_{\cC}(C_1,C_2) \to \Ext^1_{\zZ}(L_1,N)$ is surjective.
\end{LEM}
\begin{proof}
As the map in question is linear, it suffices to check that a basis for $\Ext^1_{\zZ}(L_1, N)$ is in its image.

	Let $\xi \in \Ext^1_{\zZ}(L_1, N)$ be an element such that $bac \xi$ is non-zero. First, note that by Lemma \ref{lem_theta_xi}, $\theta_{\zZ}(\xi)$ and $\theta_{\zZ}(bac \xi)$ are different elements of $\ncDef_{\zeta}(\bB, Q_{\bB})$.	As $\Psi$ is assumed to be an isomorphism, there exist $(G_1, \rho_1), (G_2, \rho_2) \in \Hom_{\wt{\Df}_3}((\bB, Q_{\bB}), (\cC, Q_{\cC}))$ such that $\Psi(G_1, \rho_1) = \theta_{\zZ}(\xi)$ and $\Psi(G_2, \rho_2) = \theta_{\zZ}(bac \xi)$. Note that short exact sequences in $\cC$ obtained by applying functors $G_1$ and $G_2$ to the short exact sequence \eqref{eqtn_res_of_S1} in $\bB$ give classes $\nu_1, \nu_2 \in \Ext^1_{\zZ}(C_1, C_2)$ such that $\theta_{\cC}(\nu_1)= (G_1, \rho_1)$ and $\theta_{\cC}(\nu_2) = (G_2, \rho_2)$. Then, as $\theta_{\zZ}T(\nu_1) = \theta_{\zZ}(\xi)$ and $\theta_{\zZ}T(\nu_2) = \theta_{\zZ}(bac\xi)$, Lemma \ref{lem_theta_xi} implies that $T(\nu_1) =(\Id + \la_1 bac)\xi$ and $T(\nu_2) = (\Id + \la_2 bac)bac \xi = bac \xi$, for some $\la_1 \in \kk$. Surjectivity of the map in question follows as $(\Id +\la_1 bac)\xi$ and $bac \xi$ span $\Ext^1_{\zZ}(L_1, N)$.
	 	\end{proof}

\vspace{0.3cm}
\subsection{A criterion of the fully faithfulness for infinitesimal extensions of functors}~\\

Let $\aA$ be an abelian category extension generated by $\aA_0 \subset \aA$
and $\Phi \colon \aA \to \zZ$ an exact functor to an abelian category. 
Since fully faithfulness of $\Phi_0  = \Phi|_{\aA_0}\colon \aA_0 \to \zZ$ implies ind-representaiblity of the deformation functor $\ncDef_{\Phi_0}$, it is worthwhile to ask when fully faithfulness of $\Phi_0$ implies that of $\Phi$.
We give a criterion for this in terms of the restriction of $\Phi$ to the first infinitesimal neighbourhood $\ngbh_1^{\aA}(\aA_0)$ of $\aA_0$ in $\aA$.

An exact functor $\Phi \colon \aA\to \zZ$ of abelian categories yields a linear map 
$$
\Phi_{A,A'}^E \colon  \Ext^1_{\aA}(A,A') \to \Ext^1_{\zZ}(\Phi_0(A), \Phi_0(A'))
$$ 
which maps $\zeta \in \Ext^1_{\aA}(A,A')$ corresponding to an exact sequence $0 \to A' \to \ol{A} \to A \to 0$ to the class of the exact sequence $0 \to \Phi(A') \to \Phi(\ol{A}) \to \Phi(A) \to 0$.

\begin{PROP}\label{prop_ff_crit}
	Consider  an abelian category $\aA$ extension generated by $\aA_0\subset \aA$. Then an exact functor $\Phi \colon \aA\to \zZ$ is fully faithful if and only if $\Phi_0:=\Phi|_{\aA_0} \colon \aA_0 \to \zZ$ is fully faithful and $\Phi_{A,A'}^E$ is a monomorphism, for any $A, A' \in \aA_0$.
\end{PROP}
\begin{proof}
	If $\Phi$ is fully faithful, so is $\Phi_0$. Further, let $0 \to A' \to A_1 \to A\to 0$ and $0 \to A' \to A_2 \to A \to 0$ be extensions which become equivalent after applying $\Phi$:
	\[
	\xymatrix{0 \ar[r] & \Phi_0(A') \ar[r] & \Phi(A_1) \ar[r] & \Phi_0(A) \ar[r] & 0 \\
		0 \ar[r] & \Phi_0(A') \ar[r] \ar[u]^{\Id} & \Phi(A_2) \ar[r] \ar[u]^{\f} & \Phi(A) \ar[r] \ar[u]^{\Id} & 0 }
	\]
	Since $\Phi$ is fully faithful, $\f = \Phi(\psi)$, for some $\psi \in \Hom(A_2, A_1)$. As $\f$ is an isomorphism, so is $\psi$, hence $A_1$ and $A_2$ are isomorphic as extensions of $A$ by $A'$. It follows that $\Phi_{A,A'}^E$ is a monomorphism.
	
	Assume now that $\Phi_0$ is fully faithful and $\Phi_{A,A'}^E$ is a monomorphism, for any $A,A'\in \aA_0$. We show by induction on $i$ that $\Phi_{A,A_i} \colon \Hom_{\aA}(A, A_i) \to \Hom_{\zZ}(\Phi(A), \Phi(A_i))$ is an isomorphism and $\Phi_{A,A_i}^E$ is a monomorphism, for any $A\in \aA_0$, $A_i \in \ngbh_i^{\aA}(\aA_0)$. 
	
	Let $0 \to A' \to A_i \to A_{i-1} \to 0$ be a short exact sequence with $A'\in \aA_0$, $A_{i-1} \in \ngbh_{i-1}^\aA(\aA_0)$.  Put $Z = \Phi(A)$, $Z'=\Phi(A')$, $Z_i= \Phi(A_i)$ and $Z_{i-1} = \Phi(A_{i-1})$. Functor 
	$\Phi$ gives a commutative diagram with exact rows:
	\begin{equation*}
		\tiny{\xymatrix{0 \ar[r]& \Hom(Z,Z') \ar[r] & \Hom(Z,Z_i) \ar[r] & \Hom(Z,Z_{i-1}) \ar[r]& \Ext^1(Z,Z') \ar[r] & \Ext^1(Z,Z_i) \ar[r] & \Ext^1(Z,Z_{i-1}) \\
			0 \ar[r]& \Hom(A,A') \ar[r] \ar[u]^{\Phi_{A,A'}} & \Hom(A,A_i) \ar[r] \ar[u]^{\Phi_{A,A_i}}& \Hom(A,A_{i-1}) \ar[r] \ar[u]^{\Phi_{A,A_{i-1}}}& \Ext^1(A,A') \ar[r] \ar[u]^{\Phi^E_{A,A'}}& \Ext^1(A,A_i) \ar[r] \ar[u]^{\Phi^E_{A,A_i}} & \Ext^1(A,A_{i-1}) \ar[u]^{\Phi^E_{A,A_{i-1}}} }}
		\end{equation*}

	By assumption, $\Phi_{A,A'}$ and $\Phi_{A,A_{i-1}}$ are isomorphisms, while $\Phi^E_{A,A'}$ and $\Phi^E_{A,A_{i-1}}$ are monomorphisms. By the five lemma, $\Phi_{A,A_i}$ is an isomorphism and $\Phi^E_{A,A_i}$ is a monomorphism.
	
	An analogous inductive argument shows that $\Phi_{A_i,A}$ is an isomorphism and $\Phi^E_{A_i,A}$ is a monomorphism, for any $A_i\in \ngbh_i^{\aA}(\aA_0)$ and $A\in \aA$.
\end{proof}

\begin{THM}\label{thm_ff_criterion}
Consider abelian categories $\aA$ and $\zZ$. Assume that $\aA$ is extension generated by a full subcategory $\aA_0 \subset \aA$. Then an exact functor $\Phi \colon \aA \to \zZ$ is fully faithful if and only if its restriction $\Phi|_{\ngbh_1(\aA_0)} \colon \ngbh_1(\aA_0) \to \zZ$ is.
\end{THM}
\begin{proof}
	Proposition \ref{prop_ff_crit} applied both to $\aA_0 \subset \aA$ and $\aA_0 \subset \ngbh_1(\aA_0)$ implies that the fully faithfulness of $\Phi$ and $\Phi|_{\ngbh_1(\aA_0)}$ is equivalent to the fully faithfulness of $\Phi|_{\aA_0}$ and injectivity of $\Phi_{A,A'}^E$, for any $A, A' \in \aA_0$. The statement follows. 
\end{proof}

\section{The null-category as the space of non-commutative deformations}\label{sec_null_cat_as_nc_def}

Let $X$ and $Y$ be normal varieties over a field $\kk$. Consider a proper morphism $f\colon X\to Y$
\begin{equation}\tag{$\ast$}
\textrm{with fibers of dimension bounded by 1, and such that }Rf_*\oO_X = \oO_Y.
\end{equation} 
Let $y \in Y$ be a closed point such that the fiber $f^{-1}(y)$ is one dimensional. Then $C := f^{-1}(y)_{\textrm{red}}$ is a 
proper algebraic curve over the field $\kk$ with $H^1(\oO_C) =0$.

To describe the structure of the reduced fiber, we introduce an \emph{incidence graph} of a reduced curve $C = \bigcup C_i$ such that 
	all of its irreducible components $C_i$ are smooth. 
	Vertices of the graph correspond to irreducible components and singular points of $C$. An edge connects a vertex corresponding to an irreducible component $C_i$ with a vertex corresponding to a singular point $c \in \textrm{Sing}(C)$ if and only if $c \in C_i$.
\begin{PROP}\cite[Theorem D.1]{BodBon}\label{prop_fiber_str} 
	Let $C$ be a reduced proper algebraic curve over a field $\kk$. Then $H^1(\oO_C) =0$ if and only if the following conditions are satisfied:
	\begin{enumerate}
			\item Every irreducible component $C_i$ of $C$ is a smooth rational curve,
			\item The incidence graph of $C$ has no cycles,
			\item The curve has normal crossing singularities.
	\end{enumerate} 
\end{PROP}

We study the \emph{null-category} of $f$: 
\begin{equation}\label{eqtn_def_A_f}
\mathscr{A}_f = \{ E\in \Coh(X)\,|\, Rf_*(E) = 0\}.
\end{equation}

We prove that $\oO_{C_i}(-1)$ are all isomorphism classes of simple objects in $\mathscr{A}_f$. We conclude that an appropriate subcategory of $\mathscr{A}_f$ ind-represents the functor of non-commutative deformations of the collection $\{\oO_{C_i}(-1)\}$. If $f$ is a morphism of threefolds we show that $\mathscr{A}_f$ in fact represents this functor.

\subsection{The null-category $\mathscr{A}_f$ under decomposition and base change}\label{ssec_A_f_under_base_change}~\\

We recall the basic properties of the null-category following \cite{BodBon}.
\begin{PROP}\cite[Proposition 2.10]{BodBon}\label{prop_g_is_exact}
	Let $f\colon X\to Y$ satisfy ($\ast$). Consider a decomposition for $f$:
	\[
	\xymatrix{X \ar[dr]^g \ar[dd]_f & \\ & Z \ar[dl]^h\\ Y&}
	\]
	Then, for $E \in \Coh(X)$ with $R^1f_* E =0$, we have $R^1g_* E = 0$. Functor $g_*$ restricts to an exact functor $g_*\colon \mathscr{A}_f \to \mathscr{A}_h$.
\end{PROP}

\begin{PROP}\cite[Proposition 2.11] {BodBon}\label{prop_vanish_Rpi_Zpi_X_of_E}
	Let $f\colon X\to Y$ satisfy $(\ast)$ and $g \colon Z \to Y$ be a morphism of schemes over a field $\kk$. Let coherent sheaf $E$ on $X$ satisfy $R^lf_* E = 0$, for $l\geq l_0$, for some $l_0 \in \{0,1\}$. Then $R^l\pi_{Z*} \pi_X^* E = 0$, for $l\geq l_0$, where $\pi_X \colon W \to X$ and $\pi_Z \colon W \to Z$ are the projections for $W = X\times_Y Z$:
	\begin{equation}\label{eqtn_diag_korean_lem}
	\xymatrix{X \times_Y Z \ar[r]^(0.6){\pi_X} \ar[d]_{\pi_Z} & X \ar[d]^f \\ Z \ar[r]^{g} & Y}
	\end{equation}
\end{PROP}

\begin{COR}\label{cor_pull-back_of_A_f}
	Let $f\colon X\to Y$ and $g \colon Z \to Y$ be as in Proposition \ref{prop_vanish_Rpi_Zpi_X_of_E}. For $E\in \mathscr{A}_f$, its pull-back $\pi_X^* E$ is an object in $\mathscr{A}_{\pi_Z}$.
\end{COR}

\vspace{0.3cm}
\subsection{Projective objects in $\mathscr{A}_f$}\label{ssec_proj_in_A_f}~\\

By \cite[Lemma 3.1]{Br1}, for $f\colon X \to Y$ satisfying $(\ast)$, the category $\mathscr{A}_f$ is the heart of the restriction of the standard \tr e on $\dD^-(\Coh(X))$ to the triangulated null-category:
$$
\cC_f =\{E^\bcdot \in \dD^-(\Coh(X))|Rf_*(E^\bcdot)  =0 \}.
$$
The semi-orthogonal decomposition $\dD^-(\Coh(X)) = \langle \cC_f, Lf^* \dD^-(\Coh(Y))\rangle$ implies existence of functor $\alpha_f^*\colon \dD^-(\Coh(X)) \to \cC_f$ left adjoint to the inclusion $\alpha_f \colon \cC_f \to \dD^-(\Coh(X))$ (see \cite[Lemma 3.1]{B}). The $\alpha_f^*\dashv \alpha_f$ adjunction unit and $Lf^*\dashv Rf_*$ adjunction counit fit into a functorial exact triangle
\begin{equation}\label{eqtn_fun_ex_tr}
Lf^*Rf_* \to \Id_{\dD^-{\Coh(X)}} \to \alpha_f \alpha_f^* \to Lf^* Rf_*[1].
\end{equation}

The standard argument shows that 
$$
\iota_f^*: = \hH^0 \circ \alpha_f^*|_{\Coh(X)} \colon \Coh(X) \to \mathscr{A}_f
$$ 
is left adjoint to the inclusion $\iota_f \colon \mathscr{A}_f \to \Coh(X)$.

To $f\colon X\to Y$ satisfying $(\ast)$ and $p\in \mathbb{Z}$, T. Bridgeland in \cite{Br1} assigned a \tr e on $\dD^b(\Coh(X))$ with  the heart $\Per{p}(X/Y)$ of $p$-perverse sheaves.
In the case when $Y$ is affine, with \cite[Proposition 3.2.5]{VdB} M. Van den Bergh constructed a vector bundle $\nN$ on $X$ which is a projective generator for the heart ${}^0 \textrm{Per}(X/Y)$ of 0-perverse \tr e on $\dD^b(\Coh(X))$. By \cite[Remark 2.6]{BodBon},
\begin{equation}\label{eqtn_def_of_P}
	\pP:=\iota_f^*\nN 
\end{equation}
is a projective generator for $\mathscr{A}_f$.

If $f\colon X\to Y$ satisfies $(\ast)$ and $Y = \Spec R$ is a spectrum of a complete Noetherian local ring, then the reduced fiber $C_{\textrm{red}} = \bigcup_{i=1}^n C_i$ of $f$ over the unique closed point $y\in Y$ is a union of rational curves $C_i$'s satisfying conditions (1)-(3) in Proposition \ref{prop_fiber_str}. In this case the Picard group of $X$ is isomorphic to $\mathbb{Z}^n$, where the isomorphism is given by the degrees of the restriction to irreducible components of $C_{\textrm{red}}$: $\lL \mapsto \deg(\lL|_{C_i})_{i=1,\ldots,n}$.

\begin{REM}(cf. \cite[Lemma 3.4.4]{VdB}\cite[Remark 2.7]{BodBon})\label{rem_divisors_D_i}
	Let $x_i \in C_i\subset X$ be a closed point such that $x_i \notin C_k$, for any $k \neq i$, and $j\colon \wt{X}_i \to \mathcal{X}_i$
	a closed embedding of the vicinity $\wt{X}_i$ of $x_i$ into a smooth scheme $\mathcal{X}_i$ over the complete local ring $R$. There exists an effective Cartier divisor $\mathcal{D}_i \subset \mathcal{X}_i$ such that scheme-theoretically
	$\mathcal{D}_i \cap j_* C_i = \{j_* x_i\}$. By pulling back $\mathcal{D}_i$ to $\wt{X}_i$, we obtain an effective divisor $D_i \subset X$ such that scheme-theoretically $D_i.C_i = \{x_i\}$ and $D_i.C_k =0$, for $k\neq i$. We denote by $\iota_{D_i} \colon D_i \to X$ the embedding of $D_i$ into $X$.

Denote by $\lL_i$ line bundles
\begin{equation}\label{eqtn_def_L_i}
\lL_i \simeq \oO_X(-D_i).
\end{equation}

Following \cite{VdB}, for every $i$, we consider vector bundle $\nN_i$:
\begin{equation}\label{eqtn_ses_def_M_i}
0 \to \lL_i \to \nN_i \to  \oO_X^{r_i-1} \to 0
\end{equation}
which corresponds to a choice of generators of $\Ext^1(\oO_X, \lL_i)$ as an $R$-module. Denote by $\nN_0=\oO_X$ the structure sheaf of $X$. By \cite[Theorem 3.5.5]{VdB}, vector bundle
$$
\nN := \bigoplus_{i=0}^n \nN_i
$$
is a projective generator for ${}^{0}\textrm{Per}(X/Y)$, i.e. $\iota_f^*\nN$ is a projective generator for $\mathscr{A}_f$.
\end{REM}

\vspace{0.3cm}
\subsection{Simple objects in $\mathscr{A}_f$}\label{ssec_simple_in_Af}~\\

\begin{PROP}\label{prop_simpl_in_A_f_supp_on_fiber}
	Let $f \colon X \to Y$ satisfy $(\ast)$ and $F \in \mathscr{A}_f$ be a simple object. Then there exists a closed point $y \in Y$ and an embedding $i\colon C_y \to X$ of the fiber of $f$ over $y$ such that $F$ is isomorphic to $i_*i^*(F)$.
\end{PROP}

\begin{proof}
	Consider a closed point $y\in Y$ such that the support of $F$ meets $C_y$. Map $\alpha \colon F \to i_*i^*(F)$ is surjective. Corollary \ref{cor_pull-back_of_A_f} implies that object $i_*i^*(F)$ is a (non-zero) object of the category $\mathscr{A}_f$. Since $F$ is simple in $\mathscr{A}_f$, morphism $\alpha$ is an isomorphism.
\end{proof}

Next we show that the category $\mathscr{A}_f$ is closed under the restriction to reduced fibers.  First, we consider a projective generator.

Note that if the fiber over a closed point of $Y$ is irreducible, then $C_\textrm{red} \simeq \mathbb{P}^1$.

\begin{LEM}\label{lem_restriction_of_P_to_red_fiber}
	Let $f\colon X\to Y$ satisfy $(\ast)$ and $Y = \Spec R$ be a spectrum of a complete Noetherian local ring. Assume that the fiber $C$ over the unique closed point of $Y$ is irreducible. Then the restriction of the projective generator $\pP$ for $\mathscr{A}_f$ to the reduced fiber $C_{\textrm{red}}$ is an object of $\mathscr{A}_f$.
\end{LEM}
\begin{proof}
	Let $D\subset X$ be a divisor as in Remark \ref{rem_divisors_D_i} and $\pP= \iota_f^*\nN$ a projective generator. 	
	Since $\alpha_f^* \oO_X \simeq 0$, applying $\iota_f^*$ to (\ref{eqtn_ses_def_M_i}) yields an isomorphism $\pP =  \iota_f^* \oO_X(-D)$. As $\iota_f^* = \hH^0 \circ \alpha_f^*$, the cohomology sheaves of the exact triangle obtained by applying (\ref{eqtn_fun_ex_tr}) to $\oO_X(-D)$ yield an exact sequence:
	$$
	f^*f_* \oO_X(-D) \to \oO_X(-D) \to \pP \to 0
	$$
	(in fact one can show that $f^*f_* \oO_X(-D) \to \oO_X(-D)$ is injective). It follows that $\pP|_{C_{\textrm{red}}}$ is a quotient of $\oO_X(-D)|_{C_{\textrm{red}}} \simeq \oO_{C_{\textrm{red}}}(-1)$. Any quotient of an invertible sheaf on a smooth curve is either an Artinian sheaf or the invertible sheaf itself (indeed, the rank of the geometric fiber of the sheaf at every closed point is bounded by 1). Since $\oO_{C_{\textrm{red}}}(-1)\in \mathscr{A}_f$ is covered by a direct sum of copies of $\pP$, the curve $C_{\textrm{red}}$ is contained in the support of $\pP$. Thus, the restriction $\pP|_{C_{\textrm{red}}}$ is not an Artinian sheaf. Hence, $\pP|_{C_{\textrm{red}}} \simeq \oO_{C_{\textrm{red}}}(-1) \in \mathscr{A}_f$, which finishes the proof.
\end{proof}

\begin{PROP}\label{prop_restr_to_red_fiber}
	Let $f\colon X \to Y$ satisfy $(\ast)$. For $F\in \mathscr{A}_f$, let $y\in Y$ be such that the fiber $C$ over $y$ is irreducible and assume that the support of $F$ meets $C$. Then the restriction of $F$ to the reduced fiber $C_{\textrm{red}}$ is an object of $\mathscr{A}_f$.
\end{PROP}
\begin{proof}
	By Corollary \ref{cor_pull-back_of_A_f}, we can assume that $F$ is supported on $C$. Therefore we can assume that $Y$ is a spectrum of a complete Noetherian local ring.
	
	The restriction map $F \to F|_{C_{\textrm{red}}}$ fits into a short exact sequence
	$$
	0 \to K_F \to F \to F|_{C_{\textrm{red}}} \to 0.
	$$
	Since morphism $f$ has fibers of relative dimension bounded by one and $R^1f_* F = 0$, we have $R^1f_* F|_{C_{\textrm{red}}} =0$. Thus, in order to prove that $F|_{C_{\textrm{red}}}\in \mathscr{A}_f$ it suffices to check that $f_* F|_{C_{\textrm{red}}} = 0$. By applying $Rf_*$ to the short exact sequence above, we get an isomorphism $f_* F|_{C_{\textrm{red}}} \simeq R^1f_* K_F$.
	
	Since the fiber $C$ is proper, $\Hom_X(\pP, F)$ is finite dimensional, for the projective generator $\pP$ for $\mathscr{A}_f$. Hence, $F$ is a quotient of a direct sum of finitely many copies of $\pP$. We shall show that $K_F$ is the quotient of finitely many copies of $K_{\pP} := \ker(\pP \to \pP|_{C_{\textrm{red}}})$, i.e. that there exists a surjective morphism $K_{\pP}^{\oplus s} \xrightarrow{\alpha} K_F$. We use Lemma \ref{lem_kernel_of_functors} below, for $\aA = \Coh(X) = \bB$, $H =  \Id_{\Coh(X)}$ and $G = (-) \otimes \oO_{C_{\textrm{red}}}$. Indeed, since $\oO_X \to \oO_{C_{\textrm{red}}}$ is an epimorphism, the morphism of functors $\Id_{\textrm{Coh}(X)} \to (-)\otimes \oO_{C_{\textrm{red}}}$ has a trivial cokernel. As both functors $\Id_{\Coh(X)}$ and $(-)\otimes \oO_{C_{\textrm{red}}}$ are right exact, assumptions of Lemma \ref{lem_kernel_of_functors} are satisfied. Thus, if $\pP^{\oplus s} \to F$ is surjective, then so is the induced map $\alpha \colon K_{\pP}^{\oplus s} \to K_F$.
	
	Since, by Lemma \ref{lem_restriction_of_P_to_red_fiber}, the sheaf $\pP|_{C_{\textrm{red}}}$ lies in $\mathscr{A}_f$, we have $R^1f_* K_{\pP} \simeq f_* \pP|_{C_\textrm{red}} =0$. As $f$ has fibers of relative dimension bounded by one, morphism $R^1f_* K_{\pP}^{\oplus s} \xrightarrow{R^1f_* \alpha} R^1f_* K_F$ is surjective. Thus, vanishing of $R^1f_* K_{\pP}$ implies that $f_* F|_{C_\textrm{red}} \simeq R^1f_* K_F \simeq 0$.
\end{proof}

Note that, the category of functors between abelian categories $\aA\to \bB$ is itself abelian. Hence, for $H, G \colon \aA \to \bB$ and morphism $\eta \colon H \to G$, there exist functors $\textrm{Ker}(\eta), \textrm{Coker}(\eta) \colon \aA \to \bB$. We have
\begin{LEM}\label{lem_kernel_of_functors}
	Let $H, G \colon \aA \to \bB$ be right exact functors and let $\eta \colon H\to G$ be a morphism with $\textrm{Coker}(\eta) \simeq 0$. Then the functor $\textrm{Ker}(\eta)$ takes surjective morphisms in $\aA$ to surjective morphisms in $\bB$.
\end{LEM}
\begin{proof}
	Let $\alpha \colon A_1 \to A_2$ be a surjective morphism in $\aA$. Morphism $\eta$ yields a morphism of short exact sequences in $\bB$:
	\[
	\xymatrix{0 \ar[r] & \textrm{Ker}(\eta)(A_2) \ar[r] & H(A_2) \ar[r]& G(A_2) \ar[r] & 0 \\
		0 \ar[r] & \textrm{Ker}(\eta)(A_1) \ar[r] \ar[u]^{\textrm{Ker}(\eta)(\alpha)} & H(A_1) \ar[r] \ar[u]^{H(\alpha)}& G(A_1) \ar[r] \ar[u]^{G(\alpha)} & 0 }
	\]
	Since $H$ and $G$ are right exact, the snake lemma yields an exact sequence
	$$
	0 \to \ker\,\textrm{Ker}(\eta)(\alpha) \to \ker\, H(\alpha) \xrightarrow{\beta} \ker \, G(\alpha) \to \textrm{coker}\, \textrm{Ker}(\eta)(\alpha) \to 0.
	$$
	It follows that $\textrm{Ker}(\eta)(\alpha)$ is surjective if and only if morphism $\beta$ above is surjective. Denote by $A'$ the kernel of $\alpha$. In diagram
	\[
	\xymatrix{G(A') \ar[rr]\ar@{->>}[dr]&  & G(A_1) \ar[r]^{G(\alpha)} & G(A_2) \ar[r] & 0 \\
		& \ker G(\alpha) \ar[ur]&&&\\
		H(A') \ar[rr] \ar[uu]^{\eta_{A'}} \ar@{->>}[dr] & & H(A_1) \ar[r]^{H(\alpha)} \ar[uu]^{\eta_{A_1}} & H(A_2) \ar[uu]^{\eta_{A_2}} \ar[r] & 0\\
		& \ker H(\alpha) \ar[ur] \ar[uu]^(0.7){\beta} &&&}
	\]
	rows are exact. Since $\eta_{A'}$ is surjective, so is morphism $\beta$, which finishes the proof.
\end{proof}

\begin{PROP}\label{prop_simpl_obj_in_A_f_for_one_contr}
	Let $f\colon X \to Y$ satisfy $(\ast)$ and assume that the fiber $C$ over a closed point of $Y$ is irreducible. Then a unique, up to isomorphism, simple object in $\mathscr{A}_f$  whose support meets $C$ is $i_* \oO_{C_{\textrm{red}}}(-1)$.
\end{PROP}
\begin{proof}
	Let $F$  be a simple object in $\mathscr{A}_f$. Proposition \ref{prop_restr_to_red_fiber} implies that $F|_{C_{\textrm{red}}}$ is an object in $\mathscr{A}_f$. Hence, the surjective morphism $F \to F|_{C_{\textrm{red}}}$ is an isomorphism. Since $Rf_* F = 0$, we have: $\hH^0(C_{\textrm{red}}, F|_{C_{\textrm{red}}}) \simeq 0 \simeq \hH^1(C_{\textrm{red}}, F|_{C_{\textrm{red}}})$. As $C_{\textrm{red}} \simeq \mathbb{P}^1$, it follows that $F|_{C_{\textrm{red}}} \simeq \oO_{C_{\textrm{red}}}(-1)$.
\end{proof}

Now we describe simple objects in $\mathscr{A}_f$, for the case when the fiber has many components. 
First we consider the case of a morphism $f\colon X\to Y$ which satisfies $(\ast)$ and $Y = \Spec R$ is the specturm of a complete Noetherian local ring. We denote by $C$ the fiber of $f$ over the unique closed point of $Y$ and by $i\colon C\to X$ the closed embedding. We consider effective divisors $\iota_{D_i} \colon D_i \to X$ as in Remark \ref{rem_divisors_D_i}.

For $E\in \Coh(X)$, by $\textrm{Tor}_0(E)$ we denote the maximal subsheaf of $E$ with a zero-dimensional support. The length of an Artinian sheaf is its dimension as $\kk$-vector space

\begin{LEM}\label{lem_rank_on_A_f}
	Let $f\colon X\to Y$ satisfy $(\ast)$ and $Y = \Spec R$ be a spectrum of a complete Noetherian local ring. Consider a coherent sheaf $F$ supported on $C$ with $\textrm{Tor}_0(F)=0$. Then $L^1\iota_{D_i}^*F =0$ and the length of $\iota_{D_i}^*F$ does not depend on the choice of point $x_i \in C_i\setminus \cup_{j\neq i} C_j$ and divisor $D_i$.
\end{LEM}

\begin{proof}
	Isomorphism ${\iota_{D_i}}_* L\iota_{D_i}^* F\simeq \oO_{D_i} \otimes^L F$ and short exact sequence
	\begin{equation*}
	0 \to \oO_{X}(-D_i) \to \oO_{X} \to \oO_{D_i} \to 0
	\end{equation*}
	imply that ${\iota_{D_i}}_*L^1\iota_{D_i}^*(F)$ is the kernel of the morphism $F(-D_i) \to F$. The kernel would have a zero dimensional support, hence it is zero by assumption on $F$.
	
	Therefore, $\iota_{D_i}^* F \simeq L\iota_{D_i}^*F$, which implies that the length of $\iota_{D_i}^*(F)$ equals the Euler characteristic $\chi(H^{\bcdot}(X, \oO_{D_i}\otimes^L F))$. Hence it depends only on the classes of $F$ and $\oO_{D_i}$ in Grothendieck group $K_0(X)$. This implies that the length is independent of the point and the divisor.
\end{proof}
Note that, for any non-zero sheaf $E'$ with zero-dimensional support, the sheaf $f_*(E')$ is non-zero. Since for any coherent sheaf $E$, the direct image $f_* \textrm{Tor}_0(E)$ is a subsheaf of $f_*(E)$, we have $f_* \textrm{Tor}_0(E)=0$, for any $E \in \mathscr{A}_f$. Hence, $\textrm{Tor}_0(E)=0$, for such $E$.

Consider the full abelian subcategory in $\mathscr{A}_f$:
\begin{equation}\label{eqtn_def_A_fC}
\mathscr{A}_{f,C} = \{E \in \mathscr{A}_{f}\,|\, \textrm{Supp}\,E \subset C\}.
\end{equation}
Lemma \ref{lem_rank_on_A_f} allows us to unambiguously define numbers
$$
r_i(F) = \textrm{length}\, \iota_{D_i}^*(F),
$$
for any $F \in \mathscr{A}_{f,C}$.

\begin{PROP}\label{prop_A_fE_has_objc_with_fin_length}
	Functions $r_i$ are well-defined on the Grothendieck group of $\mathscr{A}_f$. If $F\in \mathscr{A}_{f,C}$ and $r_i(F) = 0$, for every $i$, then $F\simeq 0$. Moreover, $\mathscr{A}_{f,C}$ is of finite length.
\end{PROP}

\begin{proof}
	Since $C\subset X$ is proper, category $\mathscr{A}_f$ is Hom- and $\Ext^1$-finite.
	
	As $\mathscr{A}_{f,C}$ does not contain sheaves with zero-dimensional support, any $F \in \mathscr{A}_{f,C}$ with $r_i(F)$ equal to zero, for all $i$, must necessarily be the zero sheaf. Since $L^1\iota_{D_i}^*(F)$ is zero for $F\in \mathscr{A}_{f,C}$ by Lemma \ref{lem_rank_on_A_f}, numbers $r_i(F)$ are additive on short exact sequences in $\mathscr{A}_{f,C}$.  This implies that the length of every object is finite. 
\end{proof}

Let $f\colon X\to Y$ satisfy $(\ast)$. For a closed point $y \in Y$ and the fiber $C_y$ over $y$ we put
$$
\mathscr{A}_{f,C_y} = \{E \in \mathscr{A}_f\,|\, \textrm{Supp }E \subset C_y\}.
$$

\begin{THM}\label{thm_simple_objects_in_A_f}
	Let $f \colon X\to Y$ satisfy $(\ast)$. For closed $y\in Y$, denote by $\iota_{y,i} \colon C_{y,i} \to X$ the embeddings of irreducible components of $C_{y,\textrm{red}}$. Then $\{{\iota_{y,i}}_* \oO_{C_{y,i}}(-1)\}_{y \in Y}$, respectively $\{{\iota_{y,i}}_* \oO_{C_{y,i}}(-1)\}$, is the set of all isomorphism classes of simple objects in $\mathscr{A}_f$, respectively in $\mathscr{A}_{f,C_y}$. 	
\end{THM}

\begin{proof}
	By Proposition \ref{prop_simpl_in_A_f_supp_on_fiber} any simple object $F\in \mathscr{A}_f$ is isomorphic to $\iota_{y*}\overline{F}$, for some closed point $y\in Y$, embedding $\iota_y \colon C_y\to X$ and object $\overline{F} \in \mathscr{A}_{f,C_y}$. Since functor $\iota_{y*}\colon \mathscr{A}_{f,C_y} \to \mathscr{A}_f$ is exact and has no kernel, object $\overline{F}$ is simple in $\mathscr{A}_{f,C_y}$. Therefore we may assume that $Y= \Spec(R)$ is a spectrum of a complete Noetherian local ring, $C$ is the fiber of $f$ over the unique closed point of $Y$, and $\iota_i \colon C_i \to X$ is the embedding of an irreducible component of $C_\textrm{red}$.
	
	Clearly, objects ${\iota_i}_*\oO_{C_i}(-1)$ are simple as any proper quotient $E$ would necessarily satisfy $r_j(E) =0$, for all $j$.

	Let $F \in \mathscr{A}_{f}$ be simple and assume that $r_i(F) \neq 0$, for some $i$. Let $Z$ be the normalization of $ \textrm{Proj}_{Y} \bigoplus_{l\geq 0} f_*(\oO_{X}(lD_i))$ and let $h \colon Z \to Y$ denote the canonical morphism. There exists a rational map $\wt{g}\colon X\to \textrm{Proj}_{Y} \bigoplus_{l\geq 0} f_*(\oO_{X}(lD_i))$ which takes a point $x \in X$ to the ideal of sections of $\bigoplus\oO_X(lD_i)$ that vanish at $x$. Since divisor $D_i$ in the linear system $|D_i|$ can be chosen in such a way that its unique closed point is any given $x_i \in C_i \setminus \bigcup_{j\neq i} C_j$, base locus of $|D_i|$ is empty and morphism $\wt{g}$ is well-defined on all closed points of $X$. Since $X$ is normal, map $\wt{g}$ admits a lift to $g\colon X \to Z$. As $lD_i.C_j = 0$, for $j \neq i$, morphism $g$ contracts all components of the fiber of $f$ but $C_i$.
	
	Since $f$ and $hg$ are birational morphisms that coincide on a dense open set, they are equal. Hence, map $f \colon X \to Y$ can be decomposed as
	\[
	\xymatrix{X \ar[dr]^{g} \ar[dd]_{f} & \\ & Z \ar[dl]^{h} \\ Y &}
	\]
	Since $g$ is proper and $Z$ is normal, $g_*\oO_X \simeq \oO_Z$ (see \cite[Lemma 4.1]{BodBon}). Fibers of $g$ and $h$ are of relative dimension bounded by one, hence $Rg_* \oO_X \simeq \oO_Z$ (see Proposition \ref{prop_g_is_exact}). It follows that $Rh_* \oO_Z \simeq Rh_* Rg_* \oO_X \simeq\oO_Y$, i.e. morphism $h$ satisfies $(\ast)$. From the decomposition of $f$ it follows that the fiber of $h$ over the closed point of $Y$ is irreducible.
	
	Since morphism $g$ takes the component $C_i$ onto the fiber of $h$, the sheaf $g_* F$ is non-zero. Propositions \ref{prop_g_is_exact} and \ref{prop_A_fE_has_objc_with_fin_length} implies that $g_* F \in \mathscr{A}_h$. It follows from Proposition \ref{prop_simpl_obj_in_A_f_for_one_contr} that there exists an injective morphism $\oO_{C_Z}(-1)\to g_*(F)$, for the reduced fiber $C_Z$ of $h$.
	
	By adjunction, there exists a non-zero morphism $\alpha :g^* \oO_{C_Z}(-1)\to F$. Sheaf $g^*\oO_{C_Z}(-1)$ is an object in $\mathscr{A}_f$. Indeed, we have: $Rf_* Lg^*g_* \oO_{C_i}(-1) = Rh_* g_* \oO_{C_i}(-1) =  0$, as $g_* \oO_{C_i}(-1) \in \mathscr{A}_h$. \cite[Lemma 2.9]{BodBon} implies that $L^ig^*g_* \oO_{C_i}(-1) \in \mathscr{A}_f$, for $i>0$. It then follows from the spectral sequence $R^qf_*L^sg^*g_* \oO_{C_i}(-1) \Rightarrow R^{q-s}f_*Lg^*g_* \oO_{C_i}(-1) =0$, that $g^*g_* \oO_{C_i}(-1)$ is an object in $\mathscr{A}_f$ too.
	
	Since $F$ is simple in $\mathscr{A}_f$, map $\alpha$ is surjective and fits into a short exact sequence
	\begin{equation}\label{eqtn_ses_for_simple}
	0 \to A_1 \to g^*\oO_{C_Z}(-1) \to F \to 0
	\end{equation}
	with $A_1 \in \mathscr{A}_{g}$. Since $r_i(g^*\oO_{C_Z}(-1)) =1$ and $F$ is a quotient of $g^*(\oO_{C_Z}(-1))$, we have $r_i(F)=1$. It follows that sheaf $A_1$ is supported on the union of the components of the fiber of $f$ different from $C_i$.
	
	In view of adjunction
	\begin{align*}
	&\Hom(g^*\oO_{C_Z}(-1), {\iota_i}_*\oO_{C_i}(-1)) \simeq \Hom(\oO_{C_Z}(-1),g_*{\iota_i}_*\oO_{C_i}(-1)) \\ &\simeq \Hom(\oO_{C_Z}(-1), \oO_{C_Z}(-1)),
	\end{align*}
	we have a non-zero morphism $\beta \colon g^*\oO_{C_Z}(-1) \to {\iota_i}_*\oO_{C_i}(-1)$. On the other hand the space $\Hom(A_1, {\iota_i}_*\oO_{C_i}(-1)) \simeq \Hom(\iota_i^*A_1, \oO_{C_i}(-1))$ is zero, because the support  of $\iota_i^*(A_1)$ is contained in $C_i \cap (\bigcup_{j \neq i} C_j)$.
	
	Hence, applying $\Hom(-, \iota_{i*} \oO_{C_i}(-1))$ to sequence (\ref{eqtn_ses_for_simple}) implies that morphism $\beta$ factors through a non-zero morphism $F \to {\iota_i}_*\oO_{C_i}(-1)$, which is necessarily an isomorphism, as both $F$ and ${\iota_i}_*\oO_{C_i}(-1)$ are simple objects in $\mathscr{A}_{f}$.
\end{proof}

Category $\dD^b(\Coh(X))$ admits also $-1$ perverse \tr e (\emph{cf.} \cite{Br1, VdB}) with heart $\Per{-1}(X/Y)$. Then
$$
\mathscr{A}_f[1]= \{E \in  \Per{-1}(X/Y)\,|\, Rf_*E =0\}.
$$
Since the functor $Rf_* \colon \Per{-1}(X/Y) \to \Coh(Y)$ is exact, $ \mathscr{A}_f[1] \subset \Per{-1}(X/Y)$ is closed under subobjects and quotient objects. Hence, for a simple object $S \in \mathscr{A}_f$, the shift $S[1]$ is simple in $\Per{-1}(X/Y)$. Using the classification of irreducible projective objects in $\Per{-1}(X/Y)$ \cite[Proposition 3.5.4]{VdB} M. Van den Bergh described simple object in $\Per{-1}(X/Y)$ as $\oO_{C_i}(-1)[1]$ and one extra object, the structure sheaf of the schematic fiber over the closed point, \cite[Proposition 3.5.7]{VdB}. Since the structure sheaf of the schematic fiber does not lie in $\mathscr{A}_f$, sheaves $\oO_{C_i}(-1)$ are all isomorphism classes of simple objects in $\mathscr{A}_f$. The above argument provides an alternative proof of Theorem \ref{thm_simple_objects_in_A_f}.

\vspace{0.3cm}
\subsection{The null-category $\mathscr{A}_{f,C}$ ind-represents deformations of $\oO_{C_i}(-1)$}\label{ssec_A_f_C_ind-rep}~\\

Let $f\colon X\to Y$ satisfy $(\ast)$. Let $C = \bigcup_{i=1}^n C_i$ be the fiber over a closed point $y\in Y$ (see Proposition \ref{prop_fiber_str}).
We consider the $\kk^{\oplus n}$-object $(\oO_{C_1}(-1),\ldots, \oO_{C_n}(-1))$ and the functor $\gamma \colon \modu \kk^{\oplus n} \to \Coh(X)$ corresponding to it.

\begin{THM}\label{thm_universal_def}
	Let $f\colon X\to Y$ satisfy $(\ast)$ and $y\in Y$ be a closed point. The 
	category $\mathscr{A}_{f,C}$ ind-represents the functor $\ncDef_{\gamma}$. 
\end{THM}
\begin{proof}
	By Proposition \ref{prop_A_fE_has_objc_with_fin_length} and Theorem \ref{thm_simple_objects_in_A_f} category $\mathscr{A}_{f,C}$ is of finite length with simple objects $\oO_{C_i}(-1)$. The statement follows from Theorem \ref{thm_pro-represent_for_simple_coll}. 
\end{proof}	
\begin{COR}\label{cor_univ_def_dim_3}
	Let $f\colon X\to Y$ be a birational morphism satisfying $(\ast)$. If $Y$ is a spectrum of a complete Noetherian local ring and the exceptional locus of $f$ is a curve $C$ contracted to a point, then the category $\mathscr{A}_f$ represents the functor $\ncDef_{\gamma}$.
\end{COR}	
\begin{proof}
	The assumptions on $f$ imply that $\mathscr{A}_f\simeq \mathscr{A}_{f,C}$, for the fiber $C$ of $f$ over the closed point $y\in Y$.
	Since  $\pP\in \mathscr{A}_f$ as in (\ref{eqtn_def_of_P}) is projective, $\mathscr{A}_f \in \Df_n$ is a Deligne finite category. We conclude by Theorem \ref{thm_universal_def}. 
\end{proof}

\appendix

\bibliographystyle{alpha}
\bibliography{../../ref}

\def\cprime{$'$} \def\cprime{$'$}
\begin{thebibliography}{VdB04}

\bibitem[ASS06]{ASS}
I.~Assem, D.~Simson, and A.~Skowro{\'n}ski.
\newblock {\em Elements of the representation theory of associative algebras.
  {V}ol. 1}, volume~65 of {\em London Mathematical Society Student Texts}.
\newblock Cambridge University Press, Cambridge, 2006.
\newblock Techniques of representation theory.

\bibitem[AZ01]{ArtZha}
M.~Artin and J.~J. Zhang.
\newblock Abstract {H}ilbert schemes.
\newblock {\em Algebr. Represent. Theory}, 4(4):305--394, 2001.

\bibitem[Bas60]{Bass}
H.~Bass.
\newblock Finitistic dimension and a homological generalization of semi-primary
  rings.
\newblock {\em Trans. Amer. Math. Soc.}, 95:466--488, 1960.

\bibitem[BB22]{BodBon}
A.~Bodzenta and A.~Bondal.
\newblock Flops and spherical functors.
\newblock {\em Compos. Math.}, 158(5):1125--1187, 2022.

\bibitem[Bon89]{B}
A.~I. Bondal.
\newblock Representations of associative algebras and coherent sheaves.
\newblock {\em Izv. Akad. Nauk SSSR Ser. Mat.}, 53(1):25--44, 1989.

\bibitem[Bri02]{Br1}
T.~Bridgeland.
\newblock Flops and derived categories.
\newblock {\em Invent. Math.}, 147(3):613--632, 2002.

\bibitem[B{\"u}h10]{Buehl}
T.~B{\"u}hler.
\newblock Exact categories.
\newblock {\em Expo. Math.}, 28(1):1--69, 2010.

\bibitem[Del90]{Del}
P.~Deligne.
\newblock Cat\'egories tannakiennes.
\newblock In {\em The {G}rothendieck {F}estschrift, {V}ol.\ {II}}, volume~87 of
  {\em Progr. Math.}, pages 111--195. Birkh\"auser Boston, Boston, MA, 1990.

\bibitem[DW16]{DonWem}
W.~Donovan and M.~Wemyss.
\newblock Noncommutative deformations and flops.
\newblock {\em Duke Math. J.}, 165(8):1397--1474, 2016.

\bibitem[DW19]{DonWem1}
W.~Donovan and M.~Wemyss.
\newblock Twists and braids for general 3-fold flops.
\newblock {\em J. Eur. Math. Soc. (JEMS)}, 21(6):1641--1701, 2019.

\bibitem[Gab62]{Gabriel}
P.~Gabriel.
\newblock Des cat\'egories ab\'eliennes.
\newblock {\em Bull. Soc. Math. France}, 90:323--448, 1962.

\bibitem[Gro61]{EGAII}
A.~Grothendieck.
\newblock \'{E}l\'ements de g\'eom\'etrie alg\'ebrique. {II}. \'{E}tude globale
  \'el\'ementaire de quelques classes de morphismes.
\newblock {\em Inst. Hautes \'Etudes Sci. Publ. Math.}, (8):222, 1961.

\bibitem[Kaw18]{Kaw5}
Y.~Kawamata.
\newblock On multi-pointed non-commutative deformations and {C}alabi-{Y}au
  threefolds.
\newblock {\em Compos. Math.}, 154(9):1815--1842, 2018.

\bibitem[Kel90]{Kel4}
B.~Keller.
\newblock Chain complexes and stable categories.
\newblock {\em Manuscripta Math.}, 67(4):379--417, 1990.

\bibitem[Koh68]{Koh}
K.~Koh.
\newblock On a structure theorem for a semi-simple ring-like object.
\newblock {\em J. Algebra}, 10:360--367, 1968.

\bibitem[KS06]{KasSch2}
M.~Kashiwara and P.~Schapira.
\newblock {\em Categories and sheaves}, volume 332 of {\em Grundlehren der
  Mathematischen Wissenschaften [Fundamental Principles of Mathematical
  Sciences]}.
\newblock Springer-Verlag, Berlin, 2006.

\bibitem[Lau02]{Lau}
O.~A. Laudal.
\newblock Noncommutative deformations of modules.
\newblock {\em Homology Homotopy Appl.}, 4(2, part 2):357--396, 2002.
\newblock The Roos Festschrift volume, 2.

\bibitem[Pop73]{Pop}
N.~Popescu.
\newblock {\em Abelian categories with applications to rings and modules}.
\newblock Academic Press, London-New York, 1973.
\newblock London Mathematical Society Monographs, No. 3.

\bibitem[Qui73]{Quillen}
Daniel Quillen.
\newblock Higher algebraic {$K$}-theory. {I}.
\newblock pages 85--147. Lecture Notes in Math., Vol. 341, 1973.

\bibitem[Sch68]{Schle}
M.~Schlessinger.
\newblock Functors of {A}rtin rings.
\newblock {\em Trans. Amer. Math. Soc.}, 130:208--222, 1968.

\bibitem[Tod07]{Toda3}
Y.~Toda.
\newblock On a certain generalization of spherical twists.
\newblock {\em Bulletin de la Soci\'{e}t\'{e} Math\'{e}matique de France},
  135:119--134, 2007.

\bibitem[VdB04]{VdB}
M.~Van~den Bergh.
\newblock Three-dimensional flops and noncommutative rings.
\newblock {\em Duke Math. J.}, 122(3):423--455, 2004.

\end{thebibliography}
\end{document}